\documentclass[a4paper, 11pt]{article}

\usepackage[titletoc,toc,title]{appendix}

\input xy   
\xyoption{all}
\usepackage{amsmath, amsthm, amsfonts, amssymb}
\usepackage{enumerate}
\usepackage{stackrel}
\usepackage{colordvi}
\usepackage{hyperref}
\usepackage[english]{babel}
\usepackage{color}
\usepackage{verbatim}
\usepackage[all]{xy}
\usepackage{graphicx}
\usepackage{mathrsfs}

\DeclareMathAlphabet{\mathcalligra}{T1}{calligra}{m}{n}
\DeclareMathAlphabet{\mathscrmin}{T1}{mathscr}{m}{n}

\numberwithin{equation}{section}
\theoremstyle{plain}
  \begingroup
        \newtheorem{theorem}[equation]{Theorem}
        \newtheorem{lemma}[equation]{Lemma}
        \newtheorem{proposition}[equation]{Proposition}
        \newtheorem{corollary}[equation]{Corollary}

\endgroup

\usepackage[textwidth=0.89in,textsize=scriptsize]{todonotes}

\usepackage{amsmath,amssymb,amsthm, xspace,a4wide,color}
\theoremstyle{definition}
  \begingroup
        \newtheorem{remark}[equation]{Remark}
        
        \newtheorem{sinnadastandard}[equation]{}

	\newtheorem{notation}[equation]{Notation}
\endgroup

\newcommand{\cc}{\mathcal}

\newcommand{\ff}{\mathsf}

\newcommand{\de}{definition}
\newcommand{\prop}{proposition}
\newcommand{\A}{\mathcal{A}}
\newcommand{\B}{\mathcal{B}}

\newcommand{\Cat}{\mathcal{C}at}

\newcommand{\mr}[1]{\overset {#1} {\longrightarrow}}

\newcommand{\xr}[1]{\xrightarrow {#1}}

\newcommand{\Mr}[1]{\overset {#1} {\Longrightarrow}}

\newcommand{\mrpairviejo}[2]
   {
    \xymatrix@C=5ex@R=2.4ex
            {
             {} \ar@<1.6ex>[r]^{#1} 
	            \ar@<-1.1ex>[r]^{#2} 
	         & {}
            }
   }

\newcommand{\mrpair}[2]
   {
    \xymatrix@C=5ex@R=2.4ex
            {
             {} \ar@<1ex>[r]^{#1} 
	            \ar@<-1ex>[r]_{#2} 
	         & {}
            }
   }

 \newcommand{\mrpairc}[2]
   {
    \xymatrix@C=5ex@R=2.4ex
            {
             {} \ar@<1ex>[r]^{#1} 
	            \ar@<-1ex>[r]|{o}_{#2} 
	         & {}
            }
   }

 \newcommand{\mrpaircc}[2]
   {
    \xymatrix@C=5ex@R=2.4ex
            {
             {} \ar@<1ex>[r]|{o}^{#1} 
	            \ar@<-1ex>[r]|{o}_{#2} 
	         & {}
            }
   }

\newcommand{\mlpair}[2]
   {
    \xymatrix@C=5ex@R=2.4ex
            {
             {} 
              & {} \ar@<1.0ex>[l]_{#2} 
	          \ar@<-1.7ex>[l]_{#1}
            }
    }

\newcommand{\cellrd}[3] 
 {
  \xymatrix@C=7ex@R=2.4ex
         {
          {} \ar@<1.6ex>[r]^{#1} 
             \ar@{}@<-1.3ex>[r]^{\!\! {#2} \, \!\Downarrow}
             \ar@<-1.1ex>[r]_{#3} 
          & {}
         }
 }
 \newcommand{\modif}[3] 
 {
  \xymatrix@C=7ex@R=2.4ex
         {
          {} \ar@<1.6ex>@{=>}[r]^{#1} 
             \ar@{}@<-1.3ex>@{=>}[r]^{\!\! {#2} \, \!\downarrow}
             \ar@<-1.1ex>[r]_{#3} 
          & {}
         }
 }
 \newcommand{\scellrd}[3] 
 {
  \xymatrix@C=4.5ex@R=2.4ex
         {
          {} \ar@<1.6ex>[r]^{#1}
             \ar@{}@<-1.3ex>[r]^{\!\! {#2} \, \!\Downarrow}
             \ar@<-1.1ex>[r]_{#3}
          & {}
         }
}

\newcommand{\cellld}[3] 
 {
  \xymatrix@C=6ex@R=2.4ex
         {
            {} 
          & {} \ar@<1.0ex>[l]^{#3} 
          \ar@{}@<-1.7ex>[l]^{\!\! {#2} \, \!\Downarrow}
	                                 \ar@<-1.7ex>[l]_{#1}
         }
 }

\newcommand{\cellpairrd}[4] 
 {
  \xymatrix@C=10ex@R=2.4ex
         {
          {} \ar@<1.6ex>[r]^{#1} 
             \ar@{}@<-1.3ex>[r]^{\!\! {#2} \, \!\Downarrow 
                                 \;\; {#3} \, \!\Downarrow }
             \ar@<-1.1ex>[r]_{#4} 
          & {}
         }
 }

\newcommand{\cellpairrdc}[4] 
 {
  \xymatrix@C=10ex@R=2.4ex
         {
          {} \ar@<1.6ex>[r]^{#1} 
             \ar@{}@<-1.3ex>[r]^{\!\! {#2} \, \!\Downarrow 
                                 \;\; {#3} \, \!\Downarrow }
             \ar@<-1.1ex>[r]|{o}_{#4} 
          & {}
         }
 }
 
\newcommand{\cellpairrdcc}[4] 
 {
  \xymatrix@C=10ex@R=2.4ex
         {
          {} \ar@<1.6ex>[r]|{o}^{#1} 
             \ar@{}@<-1.3ex>[r]^{\!\! {#2} \, \!\Downarrow 
                                 \;\; {#3} \, \!\Downarrow }
             \ar@<-1.1ex>[r]|{o}_{#4} 
          & {}
         }
 }
 
\newcommand{\coLim}[2]
   {
    \underset{#1}{\underrightarrow{\ff{\sigma Lim}}}
    \; {#2}
   }

\newcommand{\cart}   
 {
  \mathscr{C}
 }
 
\newcommand{\coLimconw}[2]
   {
    \underset{#1}{\underrightarrow{\ff{\sigma Lim^\Sigma}}}
    \; {#2}
   }

\newcommand{\dcell}[1]  
          {
					 \ar@<8pt>@{-}[d]+<-4pt,8pt> 
           \ar@<-8pt>@{-}[d]+<4pt,8pt>
           \ar@{}[d]|{#1}
          }

\newcommand{\dcellb}[1]   
          {
           \ar@<10pt>@{-}[d]+<-5pt,8pt> 
           \ar@<-10pt>@{-}[d]+<5pt,8pt>
           \ar@{}[d]|{#1}
          }

\newcommand{\deq}        
         {
          \ar@{=}[d]
         }
        
\newcommand{\dreq}       
         {
          \ar@{=}[dr]
         }

\newcommand{\dleq}       
         {
          \ar@{=}[dl]
         }

\newcommand{\dccell}[1]    
          {
           \ar@{-}[ld] 
           \ar@{-}[rd] 
           \ar@{}[d]|{#1}  
          }

\newcommand{\dcellbb}[1]   
          {
           \ar@<20pt>@{-}[d]+<-10pt,12pt> 
           \ar@<-20pt>@{-}[d]+<10pt,12pt>
           \ar@{}[d]|{#1}
          } 

\newcommand{\dl}    
          {                        
           \ar@<-2pt>@{-}[d]+<4pt,8pt>
          }

\newcommand{\dr}    
          {                        
           \ar@<2pt>@{-}[d]+<-4pt,8pt> 
          }
\newcommand{\dc}[1]    
          {                        
           \ar@{}[d]|{#1}  
          }
\newcommand{\dcr}[1]    
          {                        
           \ar@{}[dr]|{#1}  
          }

\newcommand{\uccell}[1]      
          { 
           \ar@{-}[ur] 
           \ar@{}[u]|{#1} 
           \ar@{-}[ul] 
          }

\newcommand{\uccellb}[1]     
          { 
           \ar@<-1ex>@{-}[ur] 
           \ar@{}[u]|{#1} 
           \ar@<1ex>@{-}[ul] 
          }

\newcommand{\dcellop}[1]  
          {
					 \ar@<6pt>@{-}[d]+<6pt,8pt> 
           \ar@<-6pt>@{-}[d]+<-6pt,8pt>
           \ar@{}[d]|{#1}
          }

\newcommand{\dcellopb}[1]  
          {
					 \ar@<7pt>@{-}[d]+<7pt,8pt> 
           \ar@<-7pt>@{-}[d]+<-7pt,8pt>
           \ar@{}[d]|{#1}
          }

\newcommand{\did}{\ar@2{-}[d]}

\newcommand{\op}[1]
          {
           \ar@{-}[ld] 
           \ar@{-}[rd] 
           \ar@{}[d]|{#1}  
          }

\newcommand{\cl}[1]
          { 
           \ar@{-}[ur] 
           \ar@{}[u]|{#1} 
           \ar@{-}[ul] 
          }

\newcommand{\s}{$\sigma$}

\begin{document}
 
\title{A construction of certain weak colimits and an exactness property of the 2-category of categories}

\author{Descotte M.E., Dubuc E.J., Szyld M.}
\date{\vspace{-5ex}}

\maketitle 

\begin{abstract}

Given a 2-category $\cc{A}$, a $2$-functor $\A\mr{F} \Cat$ and a distinguished 1-subcategory 
$\Sigma \subset \cc{A}$ containing all the objects, a 
\emph{$\sigma$-cone} for $F$ (with respect to $\Sigma$) 
is a lax cone 
such that the structural $2$-cells corresponding to the arrows of $\Sigma$ are invertible. 
The conical \mbox{\s-limit} is the universal (up to isomorphism) \s-cone. 
The notion of  \s-limit generalizes the well known notions of pseudo and lax limit.
We consider the fundamental notion of \mbox{\s-filtered} pair $(\cc{A}, \; \Sigma)$ which generalizes the notion of 2-filtered 2-category. We give an explicit construction of \mbox{\s-filtered} \s-colimits of categories, construction which allows computations with these colimits.  We then state and prove a basic exactness property of the 2-category of categories, namely, that \mbox{\s-filtered} \s-colimits commute with finite weighted pseudo (or bi) limits. An important corollary of this result is that a \s-filtered \s-colimit of exact  
category valued 2-functors is exact. This corollary is essential in the 2-dimensional theory of flat and pro-representable 2-functors, that we develop elsewhere.
\end{abstract}
 

\vspace{4ex}
 
\noindent {\bf \Large Introduction.} \label{intro}

\vspace{2ex}

In this paper we develop an explicit construction of certain colimits of categories, and state and prove an important exactness property of the $2$-category $\cc{C}at$ of categories, which corresponds to the commutation of filtered colimits with finite limits of sets. 

In \S \ref{sub:terminology} we recall some necessary background and fix terminology. Given a \mbox{2-category} $\cc{A}$, a $2$-functor $\A\mr{F} \Cat$ and a distinguished 1-subcategory $\Sigma \subset \cc{A}$ containing all the objects, a 
\emph{$\sigma$-cone} for $F$ (with respect to $\Sigma$) with vertex $E$ is a lax cone 
$\{FA  \mr{\theta_A} E \}_{A\in \cc{A}}\,$, 
\mbox{$\{\theta_{B} Ff \Mr{\theta_f} \theta_A\}_{A\mr{f} B \in \A}\,$,} such that $\theta_f$ is invertible for every $f$ in $\Sigma$. In Definition \ref{def:sigmacone} we recall the notion of conical \s-colimit. It is the universal (up to isomorphism, that is, in the ``pseudo'' sense) \s-cone. Conical $\sigma$-limits are special cases of \emph{cartesian quasi limits} as considered in \mbox{\cite[I,7.9.1 iii)]{GRAY}.} For a complete definition of weighted $\sigma$-limits in terms of $\sigma$-natural transformations we refer the interested reader to \cite{DDS}. Weighted \s-limits generalize both weighted lax and pseudo-limits as considered in \cite[\S 5]{K2}. 
In Definition \ref{def:sigmafiltered} we recall the fundamental notion of \mbox{\s-filtered} pair $(\cc{A}, \; \Sigma)$ introduced in \cite{DDS}. This notion relativizes to $\Sigma$ the definition of bifiltered given in \cite{KE}.

In \S \ref{sub:DSC}, Definition \ref{olvidada1} and Theorem \ref{olvidada2},
we generalize the construction of $2$-filtered pseudo-colimits of categories of \cite{DS} and develop an explicit construction of $\sigma$-filtered $\sigma$-colimits of categories which is essential in many applications. 
In \S \ref{sub:exactness} we recall the notion of \emph{finite weight}, and prove, Theorem \ref{teo:conmutan}, that $\sigma$-filtered $\sigma$-colimits commute with finite wighted pseudo (or bi) limits in $\cc{C}at$. 

\vspace{1ex}

It is clear that the notion of \emph{exact} $2$-functor should be the preservation up to equivalence of finite weighted bilimits. An important corollary of Theorem \ref{teo:conmutan} is that a \s-filtered \s-colimit of $\Cat$ valued exact 2-functors is exact. This fact is essential in the theory of \emph{flat} $2$-functors that we develop in \cite{DDS}. 

%
%
%
%


\section{Background and terminology} \label{sub:terminology}
%
%
We refer the reader to \cite{KS} for basic notions on 2-categories. Size issues are not relevant to us here, when it is not clear from the context we indicate the smallness condition if it applies. By $\Cat$ we denote the $2$-category of (small) categories, with functors as morphisms and natural transformations as 2-cells. 
 
%
%
%
In any 2-category, we use $\circ$ to denote vertical composition and juxtaposition to denote horizontal composition. We consider juxtaposition more binding than ``$\circ$'',
thus $\alpha \beta \circ \gamma $ means $(\alpha \beta ) \circ \gamma$. 
We will abuse notation by writing $f$ instead of $id_f$ for arrows $f$ when there is no risk of confusion.
%
%
For a $2$-category $\A$ and objects $A, B \in \A$, we use the notation $\A(A,B)$ to denote the category whose objects are the morphisms between $A$ and $B$ and whose arrows are the 2-cells between those morphisms.

We fix throughout this paper a pair $(\A, \; \Sigma)$, 
where $\Sigma$ is a family of arrows of a $2$-category $\A$, containing the identities and closed under composition. We note that this amounts to a $1$-subcategory of $\A$ containing all the objects. We will use a symbol $\sigma$ accompanying a concept, it is convenient to think that $\sigma$ means that the concept is to be taken ``relative to $\Sigma$''. 
Whenever possible, we'll omit $\Sigma$ from the notation.

%
%
%
%

We now give an explicit definition of the concepts of $\sigma$-cone and (conical) $\sigma$-colimit. 
For a complete definition of weighted $\sigma$-limits in terms of $\sigma$-natural transformations we refer the interested reader to \cite{DDS}. Weighted \s-limits generalize both weighted lax and pseudo-limits as considered in \cite[\S 5]{K2}, which correspond respectively to the two extreme cases where $\Sigma$ consists of only the identities or $\Sigma$ consists of all the arrows of $\cc{A}$. 
As usual, conical $\sigma$-colimits are weighted $\sigma$-colimits where the weight is constant at the category $1$. 
An important property of conical $\sigma$-colimits is that any weighted 
$\sigma$-colimit can be expressed as a conical $\sigma$-colimit, albeit with a different pair \mbox{$(\A,\; \Sigma)$ \cite{DDS}.}
\begin{\de} \label{def:sigmacone} 
{\bf \large (\s-colimit)} Let $F:\A \mr{} \B$ be a 2-functor, and $E$ an object of $\B$. 
A \emph{$\sigma$-cone} for $F$ (with respect to $\Sigma$) with vertex $E$ is a 
lax cone $\{FA  \mr{\theta_A} E\}_{A\in \cc{A}}$,
\mbox{$\{\theta_{B} Ff \Mr{\theta_f} \theta_A \}_{A\mr{f} B \in \A}$} such that $\theta_f$ is invertible for every $f$ in $\Sigma$. 
The morphisms between two $\sigma$-cones correspond to their morphisms as lax cones. This defines the category \mbox{$Cones_\sigma^{\Sigma}(F,E)$}. 

The \emph{$\sigma$-colimit} in $\B$ (with respect to $\Sigma$) of the $2$-functor \mbox{$F:\A\mr{}\B$} is the universal 
\mbox{$\sigma$-cone,} denoted 
\mbox{$\{FA \mr{\lambda_A} \coLimconw{A\in \A}{FA}\}_{A\in \A}$,} \mbox{$\{\lambda_{B} Ff \Mr{\lambda_f} \lambda_A \}_{A\mr{f} B \in \A}$} 
in the sense that for each \mbox{$E\in \B$,} pre-composition with 
$\lambda$ is an isomorphism of categories (we omit $\Sigma$ now)
\begin{equation}\label{isoplim}
\; \B(\coLim{A\in \A}{FA},E) \mr{\lambda^*} Cones_\sigma(F,E) 
\end{equation} 
\end{\de} 
Since we need them later we make explicit below the equations that lax cones and their morphisms satisfy. 


\begin{description}
 \item[LC0] For all $A\in \A$, $\hspace{8ex} {\theta}_{id_A} = id_{{\theta}_A}$
 
 \item[LC1] For all $A \mr{f} B \mr{g} C \in \A$, $\hspace{8ex} \theta_{gf}=\theta_f \circ \theta_g Ff$
$$
\xymatrix@C=8ex
          {
              FA \ar[dr]^{{\theta}_A}  \ar[d]_{Ff}^(0.6){\; {\theta}_{f} \Uparrow} 
           \\
              FB \ar[r]^{{\theta}_B} \ar[d]_{Fg}^(0.4){\; {\theta}_{g} \Uparrow}
            & E
           \\ 
              FC \ar[ur]_{{\theta}_C}
           }
\hspace{3ex} 
\xymatrix@R=6ex{\\ \txt{=}}
\hspace{3ex}
\xymatrix@C=8ex
          {
              FA \ar[dr]^{{\theta}_A}  \ar[d]_{Ff} 
           \\
              FB \ar@{}@<-0.5ex>[r]^(0.4){{\theta}_{gf} \Uparrow} \ar[d]_{Fg}
            & E
           \\ 
              FC \ar[ur]_{{\theta}_C}
           }
$$

\item[LC2] For all $A \cellrd{f}{\gamma}{g} B \in \A$, $\hspace{8ex} \theta_f = \theta_g \circ \theta_B F\gamma$
$$
\xymatrix@C=12ex
          {
              FA \ar@<1.4ex>[dr]^{{\theta}_A}
                \ar@<-1.6ex>[d]^(0.4){F\gamma}^(0.6){\; \Rightarrow}_{Ff}  
                \ar@<1.8ex>[d]^{Fg}^(0.6){\;\;\;\;\;\;\;\;\; {\theta}_g \Uparrow} 
           \\
              FB \ar[r]_{{\theta}_B} 
            & E
           }
\hspace{3ex} 
\xymatrix@R=3ex{\\ \txt{=}}
\hspace{3ex}
\xymatrix@C=8ex
          {
              FA \ar[dr]^{{\theta}_A}  \ar[d]_{Ff}^(0.6){\; {\theta}_f \Uparrow} 
           \\
              FB \ar[r]_{{\theta}_B} 
            & E
           }
$$     

\item[LCM] For all $A \mr{f} B \in \A$, $\hspace{8ex} \theta'_f \circ \varphi_B Ff= \varphi_A \circ \theta_f$
$$
\xymatrix@C=8ex
          {
               FA \ar@<-0.8ex>[dr]_(0.25){{\theta}_A\!\!\!\!\!}
                                ^(0.45){\!\!\varphi_A}^(0.6){\Uparrow}  
                \ar@<1.8ex>[dr]^{\theta'_A} 
                \ar[d]_{Ff}^(0.8){\;\;\; {\theta}_f \Uparrow} 
           \\
              FB \ar[r]_{{\theta}_B} 
            & \hspace{0ex} E
           }
\hspace{3ex}
\xymatrix@R=3ex{\\ \txt{=}}
\hspace{3ex}
\xymatrix@C=8ex
          {
              FA \ar@<0.8ex>[dr]^{\theta'_A}  \ar[d]_{Ff}^(0.5){\; \theta'_f \Uparrow} 
           \\
              FB \ar@<1ex>[r]^{\theta'_B} \ar@<-1.6ex>[r]_(0.5){{\theta}_B}^(0.5)
                                                {\varphi_B  \Uparrow}   
            & E
           }
$$
\end{description}
We note that conical $\sigma$-limits are special cases of \emph{cartesian quasi limits} considered by J. W. Gray in \cite[I,7.9.1 iii)]{GRAY}. As it is well known for the lax case, there is an op-lax naturality involved in the definition of conical colimits as weighted colimits with weight constant $1$. That is the reason why our conical $\sigma$-colimits above have the structure 2-cells in the reverse direction than Gray's cartesian quasi colimits.

In \cite[I,7.11.4 i)]{GRAY} Gray proves that {cartesian quasi colimits} in $\Cat$ exist and gives an explicit construction of them. A {dual} proof (see \cite{DDS} for details) shows that $\sigma$-colimits in $\Cat$ exist and yields the formula, for $\A \mr{F} \Cat$, $\coLim{A \in \A}{FA} = (\pi_0(\Gamma_{F})^{op})[(\cart_\Sigma)^{-1}]$. Here $\pi_0$ is the left adjoint of the inclusion $\Cat \mr{d} 2$-$\Cat$, $(\Gamma_{F})^{op}$ is the co-2-fibration associated to $F$, and $[(\cart_\Sigma)^{-1}]$ indicates that the category of fractions is taken with respect to the arrows $(f,\varphi)$ of $\Gamma_{F}$ with $f \in \Sigma$ and $\varphi$ an isomorphism (in other words the cartesian arrows over $\Sigma$). 

However, as it is the case for filtered pseudo-colimits of categories, general categories of fractions are hard to explicitly compute and work with, and a better construction is available under extra hypothesis of filteredness on the domain category, which yield a calculus of fractions \cite{SGA2}. 
A generalization for dimension $2$ is developed in \cite{DS}. See \cite{Data} where an explicit comparison is made between the construction of \cite{DS} and a construction using the calculus of $2$-fractions of \cite{P}. 

\vspace{1ex}

We now consider a fundamental notion introduced in \cite{DDS}. This notion relativizes to $\Sigma$ the definition of bifiltered given in \cite{KE}, \cite[2.6]{DS}. See also the equivalent notion of $2$-filtered \cite{DS}. We recover these notions when $\Sigma$ consists of all the arrows of $\A$.

%

\begin{\de}{\bf \large (\s-filtered $2$-categories)}
\label{def:sigmafiltered}
 We say that a pair $(\A, \; \Sigma)$ is $\sigma$-filtered, or for brevity, that $\cc{A}$ is \s-filtered, if it is non empty and the following hold (we add a circle to an arrow $\xymatrix{\cdot \ar[r]|{o} & \cdot}$ to indicate that it belongs to $\Sigma$):
  \begin{description}
 \item[\boldmath ${\sigma F0}$] Given $A, B \in \A$, there exist $E \in \A$ and morphisms $\vcenter{\xymatrix@R=0pc{A \ar[rd]|{o}^f \\
                                                                                                   & E \\
                                                                                                   B \ar[ru]|{o}_g}}$
  
  \item[\boldmath ${\sigma F1}$] Given $\xymatrix{A \ar@<1ex>[r]^f \ar@<-1ex>[r]|{o}_g & B} \in \A$, there exist a morphism $\xymatrix{B \ar[r]|{o}^h & E}$ and a 2-cell $hf \Mr{\alpha} hg$.
  If $f \in \Sigma$, we may choose $\alpha$ invertible.
  
  \item[\boldmath ${\sigma F2}$] Given $A \cellpairrdc{f}{\alpha}{\beta}{g} B \in \A$, there exists a morphism $\xymatrix{B \ar[r]|{o}^h & E}$ such that $h \alpha = h\beta $.
  \end{description}
\end{\de}

 



For a $2$-functor $\Delta \mr{F} \A$, we say that a lax cone 
%
%
$\theta$ with vertex $E$ has \emph{arrows in $\Sigma$} if the structure arrows $F(i) \mr{\theta_i} E$ are in $\Sigma$ for all $i \in \Delta$. We refer to a $2$-functor $\Delta \mr{F} \A$ as a finite $2$-diagram if $\Delta$ is a finite $2$-category. The following is proven in \cite{DDS}:

\begin{proposition} \label{prop:filtsiicone}
A pair $(\A, \; \Sigma)$ is $\sigma$-filtered if and only if every finite $2$-diagram $\Delta \mr{F} \A$ has a 
$\sigma$-cone (with respect to $F^{-1}(\Sigma) \subset \Delta$) with arrows in $\Sigma$. 
\qed
\end{proposition}

\section{A construction of $\sigma$-filtered $\sigma$-colimits in $\Cat$} \label{sub:DSC}

In this section we will generalize the construction of $2$-filtered pseudo-colimits of categories of \cite{DS}, and construct $\sigma$-filtered $\sigma$-colimits of categories. 
Let $\A \mr{F} \Cat$ a $2$-functor, and assume that $\A$ is $\sigma$-filtered. We will construct the $\sigma$-colimit of $F$ in the sense of \mbox{Definition \ref{def:sigmacone}.} 

In a similar way to \cite{DS}, we define a quasicategory $\cc{L}_{\sigma}(F)$ and compute the $\sigma$-colimit by identifying homotopical premorphisms.
We note however that the definition of homotopy of \cite{DS} doesn't generalize in a direct way to the $\sigma$-case.
We note that, in the case where $\Sigma$ consists of all the arrows of $\A$, the definition here is equivalent to the one in \cite{DS} by \cite[1.18]{DS}.

In the same spirit, we note that our proofs of the results in this section and \S \ref{sub:exactness} below are not direct generalizations of the proofs that can be found in \cite{DS}, \cite{Canevali}. Instead, we prove a general lemma 
(\ref{lema:todofinitov2}) from which all of the other results follow. This lemma and its proof are reminiscent of the usual arguing for filtered categories and filtered colimits: one should \mbox{\emph{``go further''}} in order to obtain the equations one wishes to prove.


\begin{\de}[cf. {\cite[1.5]{DS}, Quasicategory $\cc{L}(F)$}] \label{de:quasicatLF} \ \
\begin{enumerate}
 \item Objects of $\cc{L}_{\sigma}(F)$ are pairs $(x,A)$, with $A \in \A$ and $x \in FA$, as in $\cc{L}(F)$.
 \item Premorphisms $(x,A) \xr{(u,\xi,v)} (y,B)$ consist of $\xymatrix{A \ar[r]|{o}^{u} & C}$ , $\xymatrix{B \ar[r]^{v} & C}$ and 
 \mbox{$F(u)(x) \mr{\xi} F(v)(y) \in FC$.} This is the same as in $\cc{L}(F)$, but we add the requirement that the arrow $A \mr{u} C$ is in $\Sigma$.
 \item The equivalence relation ``being homotopical to'' will be defined below, as it is more clearly expressed after an abuse of notation is introduced.
\end{enumerate}
\end{\de}
%

\begin{notation}
 We omit the letter $F$ in denoting the action of $F$ on its arguments. Thus, as in \cite[1.6]{DS}, $A \cellrd{u}{\alpha}{v} B$ indicates a $2$-cell in $\A$ as well as the corresponding natural transformation $FA \cellrd{F(u)}{F(\alpha)}{F(v)} FB$ in $\Cat$. This abuse and its use are justified
  since $F \mr{x} A$ lives in $\cc{H}om_s(\A,\Cat)^{op}$ and we abuse $F \mr{x} A$ for $\A(A,-) \mr{x} F$ (recall the Yoneda lemma).  
 
 With this notation, premorphisms are written as $2$-cells $\xymatrix@R=4ex@C=4ex@ur
         {
            F \ar[r]^x  \ar@{}[dr]|{\xi \Downarrow}  \ar[d]_y 
          & A \ar[d]|{o}^u
          \\
            B  \ar[r]_v 
          & C 
         }$. This allows one to exploit the rich interplay that exists between ``real'' $2$-cells 
         of $\A$ and premorphisms or ``fake'' $2$-cells. 
\end{notation}         
         
\begin{\de} \label{def:homotopical}
Given premorphisms $(x,A) \mrpair{(u_1,\xi_1,v_1)}{(u_2,\xi_2,v_2)} (y,B)$, 
         we say that $\xi_1$ is homotopical to $\xi_2$, and write $\xi_1 \sim \xi_2$ if there are $2$-cells $(\alpha_1,\alpha_2,\beta_1,\beta_2)$, with $\alpha_i$ invertible $i=1,2$, such that
         
         \begin{equation} \label{eq:homotopia}
         \vcenter{\xymatrix@R=0.2ex@C=3ex
         {
          & A \ar[rdd]|{o}_{u_1} \ar@/^2ex/[rrrdd]|{o}^{w_A}
          \\
          && {\;\; \txt{\scriptsize{$\alpha_1 \!\Downarrow$}}} 
          \\
          F \ar[ruu]^{x} \ar[rdd]_{y}
          & {\txt{\scriptsize{$\xi_1 \!\Downarrow$}}}
          & C_1 \ar[rr]|{o}^{w_1}
          && D
          \\
          && {\txt{\;\; \scriptsize{$\beta_1 \!\Downarrow$}}}
          \\
          & B \ar[ruu]^{v_1} \ar@/_2ex/[rrruu]|{o}_{w_B}
         }}
         =
         \vcenter{\xymatrix@R=0.2ex@C=3ex
         {
          & A \ar[rdd]|{o}_{u_2} \ar@/^2ex/[rrrdd]|{o}^{w_A}
          \\
          && {\;\; \txt{\scriptsize{$\alpha_2 \!\Downarrow$}}} 
          \\
          F \ar[ruu]^{x} \ar[rdd]_{y}
          & {\txt{\scriptsize{$\xi_2 \!\Downarrow$}}}
          & C_2 \ar[rr]|{o}^{w_2}
          && D
          \\
          && {\txt{\;\; \scriptsize{$\beta_2 \!\Downarrow$}}}
          \\
          & B \ar[ruu]^{v_2} \ar@/_2ex/[rrruu]|{o}_{w_B}
         }}
         \end{equation}
         
         We refer to the $2$-cells $(\alpha_1,\alpha_2,\beta_1,\beta_2)$ as the homotopy and write $(\alpha_1,\alpha_2,\beta_1,\beta_2): \xi_1 \Rightarrow \xi_2$.
\end{\de}         

\begin{remark} \label{rem:remarka}
 Note that for any pair of ``real'' $2$-cells $\xymatrix@R=4ex@C=4ex@ur
         {
            D \ar[r]^w  \ar@{}[dr]|{\gamma_1 \Downarrow}  \ar[d]_{w'}
          & A \ar[d]|{o}^{u_1}
          \\
            B  \ar[r]_{v_1} 
          & C 
         }$, $\xymatrix@R=4ex@C=4ex@ur
         {
            D \ar[r]^w  \ar@{}[dr]|{\gamma_2 \Downarrow}  \ar[d]_{w'}
          & A \ar[d]|{o}^{u_2}
          \\
            B  \ar[r]_{v_2} 
          & C 
         }$, by considering a $\sigma$-cone with arrows in $\Sigma$ of the diagram given by them (recall Proposition \ref{prop:filtsiicone}), we have that there exist $(\alpha_1,\alpha_2,\beta_1,\beta_2)$ satisfying an equation analogous to \eqref{eq:homotopia} above, in other words ``real $2$-cells are always homotopical''.
\end{remark}

\begin{notation} [cf. {\cite[1.7]{DS}}] \label{not:composition}
 We have the following composition of two $2$-cells over a third one
 
 $$\beta \circ_\gamma \alpha = \vcenter{\xymatrix@R=1ex@C=4ex
         {
          &
           . \ar[rd]
          \\
          &
           {\txt{\scriptsize{$\alpha \Downarrow$}}}
          &
           . \ar[rd]
          \\
           . \ar[r] \ar[rdd] \ar[ruu]
          &
           . \ar[ru] \ar[rd] 
          &
           {\txt{\scriptsize{$\gamma \Downarrow$}}}
          &
           .
          \\
          &
           {\txt{\scriptsize{$\beta \Downarrow$}}}
          &
           . \ar[ru]
          \\
          &
           . \ar[ru]
         }  }   $$
 
 To compose premorphisms $\alpha=\xi$, $\beta=\zeta$ (cf. \cite[1.10 and above]{DS}) we pick $\gamma$ from a $\sigma$-cone with arrows in $\Sigma$ of the diagram $\vcenter{\xymatrix@R=0ex@C=4ex{ & \cdot \\ B \ar[ru] \ar[rd]|{o} \\ & \cdot }}$ 
 to determine a premorphism

 \begin{equation} \label{eq:compdepremorf}
 \zeta \circ_\gamma \xi = \vcenter{\xymatrix@R=1ex@C=4ex
         {
          &
           . \ar[rd]|{o}
          \\
          &
           {\txt{\scriptsize{$\xi \Downarrow$}}}
          &
           . \ar[rd]|{o}
          \\
           F \ar[r] \ar[rdd] \ar[ruu]
          &
           B \ar[ru] \ar[rd]|{o} 
          &
           {\txt{\scriptsize{$\gamma \Downarrow$}}}
          &
           .
          \\
          &
           {\txt{\scriptsize{$\zeta \Downarrow$}}}
          &
           . \ar[ru]|{o}
          \\
          &
           . \ar[ru]
         }  }
\end{equation}
         
Given composable premorphisms $\xi_i$, $i=1,2,3$, the equality $\xi_3 \circ_\delta (\xi_2 \circ_{\gamma_1} \xi_1) = (\xi_3 \circ_{\gamma_2} \xi_2) \circ_{\eta} \xi_1$ holds if $\delta$ and $\eta$ are chosen from the diagram

\begin{equation} \label{eq:paraasociat}
 \vcenter{\xymatrix@R=1ex@C=4ex
         { & . \ar[rd] \\
         & {\txt{\scriptsize{$\xi_1 \Downarrow$}}} & . \ar[rd] \\
          & . \ar[ru] \ar[rd] & {\txt{\scriptsize{$\gamma_1 \Downarrow$}}} & . \ar[rd] \\
         F \ar[ru] \ar[rd] \ar[ruuu] \ar[rddd] & {\txt{\scriptsize{$\xi_2 \Downarrow$}}} & . \ar[ru] \ar[rd] & {\txt{\scriptsize{$\gamma_3 \Downarrow$}}} & . \\ 
          & . \ar[ru] \ar[rd] & {\txt{\scriptsize{$\gamma_2 \Downarrow$}}} & . \ar[ru] \\
         & {\txt{\scriptsize{$\xi_3 \Downarrow$}}} & . \ar[ru] \\ 
          & . \ar[ru] }}
\end{equation}

\end{notation}

The following lemma is the key to almost all the results of this paper.

\begin{lemma} \label{lema:todofinitov2}
Consider a finite sequence of functors $\A \mr{F_i} \Cat$ (note that we may have $F_i = F_j$ for $i \neq j$), premorphisms $\xymatrix@ur@R=3.5ex@C=3.5ex
         {
        F_i \ar[r]^{x_i}  \ar@{}[dr]|{\xi_i \Downarrow}  \ar[d]_{y_i} 
          & A_i \ar[d]|{o}^{u_i}
          \\
            B_i  \ar[r]_{v_i} 
          & C_i 
         }$ in $\cc{L}_{\sigma}(F_i)$
%
         and a finite set $S = \{\eta,\xi,...\}$ of finite compositions of the premorphisms over $2$-cells of $\A$ (we ask that, for any of the given compositions $\eta$, all the involved premorphisms have the same $F_i$, and all the compositions fit as in \eqref{eq:compdepremorf})
         
         Assume that a finite number of homotopy equations between these compositions hold (we ask that, for any given equation, the premorphisms on both sides of the equations have the same domain and codomain). 
         
         Then there exist $E \in \A$, and for each $i$ 
         arrows $w_i, t_i, z_i$ in $\Sigma$ 
         and $2$-cells $\mu_i$, $\nu_i$, with $\mu_i$ invertible, fitting as follows
         
         \begin{equation} \label{eq:premorfismostilde}
\widetilde{\xi_i}:
\vcenter{\xymatrix@R=0.2ex@C=3ex
         {
          & A_i \ar[rdd]|{o}_{u_i} \ar@/^2ex/[rrrdd]|{o}^{w_i}
          \\
          && {\;\; \txt{\scriptsize{$\mu_i \!\Downarrow$}}} 
          \\
          F_i \ar[ruu]^{x_i} \ar[rdd]_{y_i}
          & {\txt{\scriptsize{$\xi_i \!\Downarrow$}}}
          & C_i \ar[rr]|{o}^{z_i}
          && E,
          \\
          && {\txt{\;\; \scriptsize{$\nu_i \!\Downarrow$}}}
          \\
          & B_i \ar[ruu]^{v_i} \ar@/_2ex/[rrruu]|{o}_{t_i}
         }}
\end{equation}

\noindent such that if we consider the corresponding
premorphisms $\widetilde{\xi_i}$,
then the equations for the compositions of the $\widetilde{\xi_i}$ hold in $FE$. We mean the equations obtained by using composition in $FE$ instead of composition over the $2$-cells, and equality instead of homotopy.

Moreover, 
\begin{enumerate}
 \item The arrows $w_i$ and $t_i$ can be chosen depending only on the objects $A_i$, $B_i$ (this means precisely that if $A_i = B_j$ then $w_i = t_j$ and similarly if $A_i = A_j$ or $B_i = B_j$)
 \item The arrow $z_i$ and the $2$-cells $\mu_i$, $\nu_i$ can be chosen so that:
 \begin{enumerate}
  \item[i)] If $A_i = B_i = C_i$ and $u_i = v_i = id$, then $z_i = w_i = t_i$ and $\mu_i = \nu_i = id$.
  \item[ii)] Otherwise, i.e. considering only those pairs $(u_i,v_i)$ such that $(u_i,v_i) \neq (id,id)$, the $2$-cells $\mu_i$ and $\nu_i$ can be chosen depending only on the arrows $u_i$, $v_i$ (this means precisely that, for each  $i,j$ such that $(u_i,v_i) \neq (id,id)$ and $(u_j,v_j) \neq (id,id)$, if $u_i = v_j$ then $\mu_i = \nu_j$ and similarly if $u_i = u_j$ or $v_i = v_j$).
 \end{enumerate}
 Note that in particular, for every $i,j$, $(u_i,v_i) = (u_j,v_j)$ implies that $(\mu_i,\nu_i) = (\mu_j,\nu_j)$.
\end{enumerate}
\end{lemma}

\begin{proof}
The idea of the proof is to construct a finite $2$-diagram $\triangle \mr{G} \A$ that holds all the information of the premorphisms, the compositions over the $2$-cells and the homotopies. 
We will construct $\triangle$ (and $G$) in steps, starting with an empty category and adding objects, arrows and $2$-cells. We omit to explicitly add all the identities and all the compositions of the arrows and the $2$-cells. We note however that this procedure yields a finite $2$-category $\triangle$ because neither arrows between the same object (loops) nor $2$-cells between the same arrow ($2$-loops) can arise as compositions of the ``added'' arrows and $2$-cells.

\begin{enumerate}
 
 \item[Step 0] We begin constructing $\triangle$ by adding one object for each element of the set $\{A_i, B_i\}$ 
 consisting of all the $A_i$ and all the $B_i$. 
 Then there is no harm in labeling these objects by $A_i$, $B_i$. The reader should note that, if $A_i = A_j$, or $A_i = B_j$, or $B_i = B_j$, then the same equality holds in $\triangle$.
 \item[Step 1] We now consider all the premorphisms $\xi_i$ such that $(u_i,v_i) \neq (id,id)$. We add one new object $\star$ for each element of the set consisting of all the $C_i$ in these premorphisms, one new arrow for each element of the set $\{u_i, v_i\}$ consisting of all the arrows in these premorphisms, and define $G$ as follows:
 
 $$\vcenter{\xymatrix@R=1ex@C=4.5ex{ A_i \ar[rd]^{a} \\ & \star \\ B_i \ar[ru]_{b} }} \quad \stackrel{G}{\longmapsto} \quad \vcenter{\xymatrix@R=1ex@C=4.5ex{ A_i \ar[rd]|{o}^{u_i} \\ & C_i \\ B_i \ar[ru]_{v_i} }}$$
 
 We note that, for every premorphism $\xi_i$, there are arrows of $\triangle$ mapped to $u_i$ and $v_i$. We now want to extend this result to all the given compositions of the $\xi_i$ over the $2$-cells.
 
 \item[Step 2] 
 To simplify the notation, it is convenient to add to the set $S$ of premorphisms $\{\eta,\xi,...\}$ all the $\xi_i$ and all the compositions between them appearing in the given premorphisms $\eta, \xi,...$. For example, if $\eta = (\xi_j \circ_{\gamma} \xi_i) \circ_{\delta} \xi_m \in S$, then we add $\xi_j \circ_{\gamma} \xi_i$ to the family. When we use the letters $\eta, \xi$ we will be referring to premorphisms in this family.
 
 We will modify $\triangle$ so that
 
 \begin{sinnadastandard} \label{sin:loqueobtengo}
For each premorphism $\xymatrix@ur@R=3.5ex@C=3.5ex
         {
        F \ar[r]^{x_i}  \ar@{}[dr]|{\eta \Downarrow}  \ar[d]_{y_j} 
          & A_i \ar[d]|{o}^{u}
          \\
            B_j  \ar[r]_{v} 
          & C 
         }$ we have chosen arrows of $\triangle$, $A_i \mr{a} \star$, $B_j \mr{b} \star$, that are mapped to $u$ and $v$ respectively. We do as follows:
 \end{sinnadastandard}
 
 We modify $\triangle$, considering one $\eta$ at a time. We work inductively in the order that the composition is computed. 
 By the Step 1, \ref{sin:loqueobtengo} holds for the $\xi_i$.
  Consider then a composition 
 $$ \eta \circ_{\gamma} \xi =  
   \vcenter
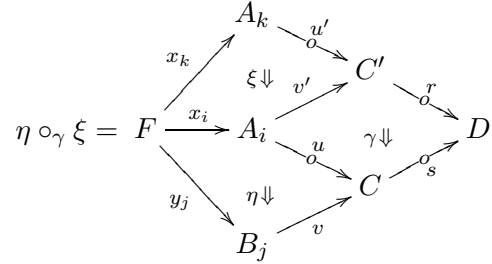
{\xymatrix@R=1ex@C=4.5ex{ & A_k \ar[rd]|{o}^{u'}  \\
 & {\;\; \txt{\scriptsize{$\xi \!\Downarrow$}}} & C' \ar[rd]|{o}^{r}  \\
 F \ar[ruu]^{x_k} \ar[r]^{x_i} \ar[rdd]_{y_j} & A_i \ar[ru]^{v'} \ar[rd]|{o}^{u} & {\;\; \txt{\scriptsize{$\gamma \!\Downarrow$}}} & D \\
 & {\;\; \txt{\scriptsize{$\eta \!\Downarrow$}}} & C \ar[ru]|{o}_{s} \\
 & B_j \ar[ru]_{v} }} 
 $$
%
%
    Then by the inductive hypothesis we have the diagram 
%
 $$\vcenter
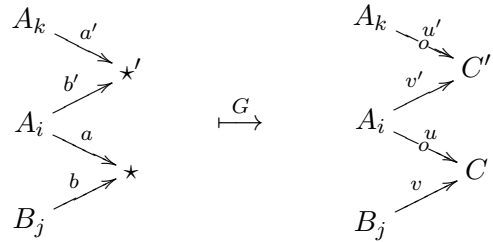
{\xymatrix@R=1ex@C=4.5ex{ & A_k \ar[rd]^{a'}  \\
 & & \star' \\
 & A_i \ar[ru]^{b'} \ar[rd]^{a} & \\
 & & \star \\
 & B_j \ar[ru]^{b} }}  
     \ \ \ \ \ \   \stackrel{G}{\longmapsto}
    \vcenter{\xymatrix@R=1ex@C=4.5ex{ & A_k \ar[rd]|{o}^{u'}  \\
 &  & C' \\
 & A_i \ar[ru]^{v'} \ar[rd]|{o}^{u} \\
 & & C \\
 & B_j \ar[ru]^{v} }} 
 $$
 We add one new object $\widetilde{\star}$, two new arrows $c,d$ and one new $2$-cell to $\triangle$, and extend $G$ as follows:
 $$\vcenter
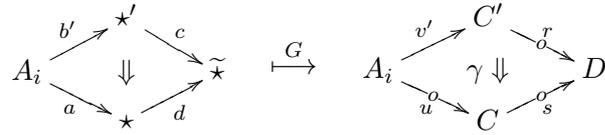
{\xymatrix@R=1ex@C=4.5ex{ & \star' \ar[rd]^c \\ A_i \ar[ru]^{b'} \ar[rd]_{a} & \Downarrow & \widetilde{\star} \\ & \star \ar[ru]_d }}  \quad \stackrel{G}{\longmapsto} \quad 
 \vcenter{\xymatrix@R=1ex@C=4.5ex{ & C' \ar[rd]|{o}^{r}  \\ A_i \ar[rd]|{o}_{u} \ar[ru]^{v'} & \gamma \Downarrow & D  \\ & C \ar[ru]|{o}_{s} }}$$
 
 Then the arrows $ca'$ and $db$ are mapped to $r u'$ and $s v$ as desired.
 
 
 We note that we add one object, two arrows and a $2$-cell once for each appearing composition over a $2$-cell $\gamma$, and not once for each $\gamma$. 
 We have to do this because two premorphisms 
 $\xi_i: \xymatrix@ur@R=1.5ex@C=1.5ex
         {
          F \ar[r] \ar@{}[dr]|{\Downarrow}
	                                                    \ar[d]    
          & \cdot \ar[d]^{u_i}
          \\
            \cdot  \ar[r]_{v_i} 
          & \cdot
         }$,
 $\xi_j: \xymatrix@ur@R=1.5ex@C=1.5ex
         {
          F \ar[r] \ar@{}[dr]|{\Downarrow}
	                                                    \ar[d]    
          & \cdot \ar[d]^{u_j}
          \\
            \cdot  \ar[r]_{v_j}
          & \cdot
         }$
 can have the same arrow $v_i = v_j = id$ but if $u_i \neq id = u_j$ then we have chosen two different arrows of $\triangle$ that map to $v_i$ and $v_j$ (cf. item 2.ii) in the statement of the lemma). Then the same $2$-cell $\gamma$ may appear in 
 two compositions $\xi \circ_{\gamma} \xi_j$ and $\xi \circ_{\gamma} \xi_i$ but must be treated differently in each case.
 \item[Step 3] We finally consider each of the homotopy equations $\eta \sim \xi$. 
 By definition, there are $2$-cells $(\alpha_1,\alpha_2,\beta_1,\beta_2)$, with $\alpha_i$ invertible $i=1,2$, such that
         \begin{equation} \label{eq:homotentreetayxi}
         \vcenter{\xymatrix@R=0.2ex@C=3ex
         {
          & A_i \ar[rdd]|{o}_{u} \ar@/^2ex/[rrrdd]|{o}^{w_{A_i}}
          \\
          && {\;\; \txt{\scriptsize{$\alpha_1 \!\Downarrow$}}} 
          \\
          F \ar[ruu]^{x_i} \ar[rdd]_{y_j}
          & {\txt{\scriptsize{$\eta \!\Downarrow$}}}
          & C \ar[rr]|{o}^{w_1}
          && D
          \\
          && {\txt{\;\; \scriptsize{$\beta_1 \!\Downarrow$}}}
          \\
          & B_j \ar[ruu]^{v} \ar@/_2ex/[rrruu]|{o}_{w_{B_j}}
         }}
         =
         \vcenter{\xymatrix@R=0.2ex@C=3ex
         {
          & A_i \ar[rdd]|{o}_{u'} \ar@/^2ex/[rrrdd]|{o}^{w_{A_i}}
          \\
          && {\;\; \txt{\scriptsize{$\alpha_2 \!\Downarrow$}}} 
          \\
          F \ar[ruu]^{x_i} \ar[rdd]_{y_j}
          & {\txt{\scriptsize{$\xi \!\Downarrow$}}}
          & C' \ar[rr]|{o}^{w_2}
          && D
          \\
          && {\txt{\;\; \scriptsize{$\beta_2 \!\Downarrow$}}}
          \\
          & B_j \ar[ruu]^{v'} \ar@/_2ex/[rrruu]|{o}_{w_{B_j}}
         }}
         \end{equation}
 	We know by \ref{sin:loqueobtengo} that there are arrows $a,b,a',b'$ of $\triangle$ mapped to $u$,$v$,$u'$,$v'$ respectively. We add one object $\ast$, four arrows $c,d,e,f$ and four $2$-cells to $\triangle$, and extend $G$ as follows:
	
	 $$ \vcenter
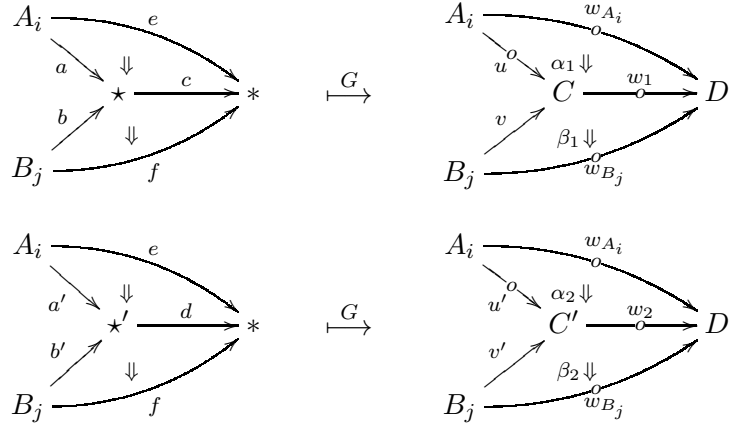
{\xymatrix@R=0.2ex@C=3ex{
          & A_i \ar[rdd]_{a} \ar@/^2ex/[rrrdd]^{e}
          \\
          && {\;\; \txt{\scriptsize{$\Downarrow$}}} 
          \\
          && \star \ar[rr]^{c}
          && \ast
          \\
          && {\txt{\;\; \scriptsize{$\Downarrow$}}}
          \\
          & B_j \ar[ruu]^{b} \ar@/_2ex/[rrruu]_{f}
         }}
	\ \ \ \ \ \stackrel{G}{\longmapsto}
	 \vcenter{\xymatrix@R=0.2ex@C=3ex
         {
          & A_i \ar[rdd]|{o}_{u} \ar@/^2ex/[rrrdd]|{o}^{w_{A_i}}
          \\
          && {\;\; \txt{\scriptsize{$\alpha_1 \!\Downarrow$}}} 
          \\
          & & C \ar[rr]|{o}^{w_1}
          && D
          \\
          && {\txt{\;\; \scriptsize{$\beta_1 \!\Downarrow$}}}
          \\
          & B_j \ar[ruu]^{v} \ar@/_2ex/[rrruu]|{o}_{w_{B_j}}
         }}$$
 
	  $$ \vcenter{\xymatrix@R=0.2ex@C=3ex{
          & A_i \ar[rdd]_{a'} \ar@/^2ex/[rrrdd]^{e}
          \\
          && {\;\; \txt{\scriptsize{$\Downarrow$}}} 
          \\
          && \star' \ar[rr]^{d}
          && \ast
          \\
          && {\txt{\;\; \scriptsize{$\Downarrow$}}}
          \\
          & B_j \ar[ruu]^{b'} \ar@/_2ex/[rrruu]_{f}
         }}
	 \ \ \ \ \ \stackrel{G}{\longmapsto}
	 \vcenter{\xymatrix@R=0.2ex@C=3ex
         {
          & A_i \ar[rdd]|{o}_{u'} \ar@/^2ex/[rrrdd]|{o}^{w_{A_i}}
          \\
          && {\;\; \txt{\scriptsize{$\alpha_2 \!\Downarrow$}}} 
          \\
          & & C' \ar[rr]|{o}^{w_2}
          && D
          \\
          && {\txt{\;\; \scriptsize{$\beta_2 \!\Downarrow$}}}
          \\
          & B_j \ar[ruu]^{v'} \ar@/_2ex/[rrruu]|{o}_{w_{B_j}}
         }}$$
\end{enumerate}

We have finished constructing $\triangle$ and $G$. By Proposition \ref{prop:filtsiicone} 
we have a $\sigma$-cone (with respect to $G^{-1}(\Sigma)$)
$\theta$ with arrows in $\Sigma$ and vertex $E$, of the diagram $\triangle \mr{G} \A$. 

Using \ref{sin:loqueobtengo} we define, for each $\eta$ 

         
 \begin{equation} \label{eq:premorfismostildeconeta}
\widetilde{\eta}:
\vcenter{\xymatrix@R=0.2ex@C=3ex
         {
          & A_i \ar[rdd]|{o}_{u} \ar@/^2ex/[rrrdd]|{o}^{\theta_{A_i}}
          \\
          && {\;\; \txt{\scriptsize{$\theta_a^{-1} \!\Downarrow$}}} 
          \\
          F \ar[ruu]^{x} \ar[rdd]_{y}
          & {\txt{\scriptsize{$\eta \!\Downarrow$}}}
          & C \ar[rr]|{o}^{\theta_{\star}}
          && E
          \\
          && {\txt{\;\; \scriptsize{$\theta_b \!\Downarrow$}}}
          \\
          & B_j \ar[ruu]^{v} \ar@/_2ex/[rrruu]|{o}_{\theta_{B_j}}
         }}
\end{equation}        
         
Note that in particular we have for each $\xi_i$, the premorphism $\widetilde{\xi_i}$ as in \eqref{eq:premorfismostilde} with 
$w_i = \theta_{A_i}$, $t_i = \theta_{B_i}$, $z_i = \theta_{\star}$, $\mu_i = \theta_{a}^{-1}$, $\nu_i = \theta_{b}$; where the object $\star$ and the arrows $a,b$ were added to $\triangle$ in step 1. By the construction of $\triangle$, items 1. and 2. at the end of the lemma hold.

We will show the following two basic properties of the $\widetilde{\eta}$ from which it will follow that the original homotopy equations hold for the $\widetilde{\xi_i}$ in $FE$, as desired. Consider two of these premorphisms $\eta$, $\xi$; then 

 
 \vspace{1ex}
 
 A. If the equation $\eta \sim \xi$ was one of the homotopy equations, then $\widetilde{\eta} = \widetilde{\xi}$.

 B. If the composition $\eta \circ_{\gamma} \xi$ appears in one of the given compositions, then $\widetilde{\eta \circ_{\gamma} \xi} = \widetilde{\eta} \circ \widetilde{\xi}$.

\vspace{1ex}

Using first A and then B, it is easy to show that the homotopy equations between the compositions of the $\xi_i$ yield the desired equations in $FE$ between the corresponding compositions of the $\widetilde{\xi_i}$. Then it only remains to prove A and B.
%
%

\vspace{1ex}

\noindent Proof of A. The homotopy $(\alpha_1,\alpha_2,\beta_1,\beta_2): \eta \Rightarrow \xi$ 
was considered in Step 3 of our construction of $\triangle$. We refer there for the notation. 
We claim that 

\begin{equation}\label{eq:widetildeetaesotracosa}
 \widetilde{\eta} =   \vcenter{\xymatrix@R=0.2ex@C=3ex
         {
          & A_i \ar[rdd]|{o}_{u} \ar@/^2ex/[rrrdd]|{o}^{w_{A_i}} \ar@/^5ex/[rrrrrdd]|>>>>>{o}^>>>>>{\theta_{A_i}}_>>>>>>>>{\theta_{e}^{-1} \Downarrow}
          \\
          && {\;\; \txt{\scriptsize{$\alpha_1 \!\Downarrow$}}} 
          \\
          F \ar[ruu]^{x_i} \ar[rdd]_{y_j}
          & {\txt{\scriptsize{$\eta \!\Downarrow$}}}
          & C \ar[rr]|{o}^{w_1}
          && D \ar[rr]|{o}^{\theta_{\ast}} && E
          \\
          && {\txt{\;\; \scriptsize{$\beta_1 \!\Downarrow$}}}
          \\
          & B_j \ar[ruu]^{v} \ar@/_2ex/[rrruu]|{o}_{w_{B_j}} \ar@/_5ex/[rrrrruu]|>>>>>{o}_>>>>>{\theta_{B_j}}^>>>>>>>>{\theta_{f} \Downarrow}
         }}
\end{equation}

	and similarly for $\widetilde{\xi}$.
	From these equalities it is clear how to show $\widetilde{\eta} = \widetilde{\xi}$ using \eqref{eq:homotentreetayxi}. So it suffices to prove \eqref{eq:widetildeetaesotracosa}.  
%
%
          By axioms {\bf LC1} and {\bf LC2} (see Definition \ref{def:sigmacone}) we have
          
          $$ \vcenter{\xymatrix@R=3ex@C=3ex{A_i \ar[rd]|{o}^u \ar@/_2ex/[rdd]|{o}_{w_{A_i}} \ar@/^5ex/[rrrrdd]|{o}^{\theta_{A_i}}_{\theta_{a} \Uparrow} \\
          \ar@{}[r]|{\stackrel{\alpha_1}{\Rightarrow}} & C \ar[d]|<<<{o}^{w_1} \ar@/^2ex/[rrrd]|{o}^{\theta_{\star}} && \\
          & D  \ar[rrr]|{o}_{\theta_{\ast}} & \ar@{}[ru]|{\theta_{c} \Uparrow} && E}}
          =
          \vcenter{\xymatrix@R=2.5ex@C=3ex{A_i \ar@/^3ex/[rrrdd]|{o}^{\theta_{A_i}} \ar[rdd]|{o}_{w_{A_i}} \\
          & \ar@{}[r]|{\theta_{e} \Uparrow} & \\
          & D  \ar[rr]|{o}_{\theta_{\ast}} && E}}$$


	    Composing with $\theta_{e}^{-1}$ and $\theta_{a}^{-1}$ it follows
	    
	  $$ \vcenter{\xymatrix@R=0.2ex@C=3ex
         {
          & A_i \ar[rdd]|{o}_{u} \ar@/^2ex/[rrrdd]|{o}^{w_{A_i}} \ar@/^5ex/[rrrrrdd]|{o}^{\theta_{A_i}}_>>>>>>>>>>{\theta_{e}^{-1} \Downarrow}
          \\
          && {\;\; \txt{\scriptsize{$\alpha_1 \!\Downarrow$}}} 
          \\
          & & C \ar[rr]|{o}^{w_1} \ar@/_4ex/[rrrr]|{o}_{\theta_{\star}}^{\theta_{c} \Downarrow}
          && D  \ar[rr]|{o}^{\theta_{\ast}} && E }} 
          =
          \vcenter{\xymatrix@R=0.2ex@C=5ex{
           A_i \ar[rd]|{o}_{u} \ar@/^3ex/[rrrd]|{o}^{\theta_{A_i}}_<<<<<<<<<<{\theta_{a}^{-1} \Downarrow} \\
           & C \ar[rr]|{o}_{\theta_{\star}} &&  E }} 
          $$

	 From this equation and the corresponding one for $\beta_1$ it is immediate to show that the diagrams in equations \eqref{eq:premorfismostildeconeta} and \eqref{eq:widetildeetaesotracosa} are equal.
           
\vspace{1ex}

\noindent Proof of B. The composition $\eta \circ_{\gamma} \xi$ was considered in Step 2 of the construction of $\triangle$, we refer there for the notation.
%
%
%
%
%
 We have to show that
  
 $$ \vcenter{\xymatrix@R=1ex@C=4.5ex{ & A_k \ar[rd]|{o}^{u'} \ar@/^4.5ex/[rrrdd]|{o}^{\theta_{A_k}} && \\
 & {\;\; \txt{\scriptsize{$\xi \!\Downarrow$}}} & C' \ar[rd]|{o}^{r}  \ar@{}[ru]|{\theta_{ca'}^{-1} \Downarrow} \\
 F \ar[ruu]^{x_k} \ar[r]^{x_i} \ar[rdd]_{y_j} & A_i \ar[ru]^{v'} \ar[rd]|{o}^{u} & {\;\; \txt{\scriptsize{$\gamma \!\Downarrow$}}} & D \ar[r]|{o}^{\theta_{\widetilde{\star}}} & E \\
 & {\;\; \txt{\scriptsize{$\eta \!\Downarrow$}}} & C \ar[ru]|{o}_{s} \ar@{}[rd]|{\theta_{db} \Downarrow} \\
 & B_j \ar[ru]_{v} \ar@/_4.5ex/[rrruu]|{o}_{\theta_{B_j}} && }}
 =
  \vcenter{\xymatrix@R=2ex@C=4.5ex{ & A_k \ar[rd]|{o}^{u'} \ar@/^4.5ex/[rrrdd]|{o}^{\theta_{A_k}} && \\
 & {\;\; \txt{\scriptsize{$\xi \!\Downarrow$}}} & C' \ar[rrd]|{o}^{\theta_{\star'}}    \ar@{}[ru]|{\theta_{a'}^{-1} \Downarrow \;\;} \\
 F \ar[ruu]^{x_k} \ar[r]^{x_i} \ar[rdd]_{y_j} & A_i \ar[ru]^{v'} \ar[rd]|{o}_{u} \ar[rrr]|{o}|>>>>>>>>{\theta_{A_i}}^{\theta_{b'} \Downarrow}_{\theta_{a}^{-1} \Downarrow}   &&& E \\
 & {\;\; \txt{\scriptsize{$\eta \!\Downarrow$}}} & C \ar[rru]|{o}_{\theta_{\star}} \ar@{}[rd]|{\theta_{b} \Downarrow} \\
 & B_j \ar[ru]_{v} \ar@/_4.5ex/[rrruu]|{o}_{\theta_{B_j}} && }} $$
 
 By axiom {\bf LC1} (see Definition \ref{def:sigmacone}) it suffices to show 

 \begin{equation} \label{eq:toshow}
  \vcenter{\xymatrix@R=2.5ex@C=4.5ex{ 
  & C' \ar[rd]|{o}^{r} \ar@/^3ex/[rrrd]|{o}^{\theta_{\star'}} \ar@{}[rrrd]|{\theta_{c}^{-1} \Downarrow}  \\
  A_i \ar[ru]^{v'} \ar[rd]|{o}_{u} & {\;\; \txt{\scriptsize{$\gamma \!\Downarrow$}}} & D \ar[rr]|{o}^{\theta_{\widetilde{\star}}} && E \\
  & C \ar[ru]|{o}_{s} \ar@/_3ex/[rrru]|{o}_{\theta_{\star}} \ar@{}[rrru]|{\theta_{d} \Downarrow} }}  
    =
    \vcenter{\xymatrix@R=2.5ex@C=4.5ex{
  & C' \ar@/^2ex/[rrd]|{o}^{\theta_{\star'}} && \\
  A_i \ar[ru]^{v'} \ar[rd]|{o}_{u} \ar[rrr]|{o}|>>>>>>>>{\theta_{A_i}} \ar@{}[rrru]|{\theta_{b'} \Downarrow} \ar@{}[rrrd]|{\theta_{a}^{-1} \Downarrow}   &&& E \\
   & C \ar@/_2ex/[rru]|{o}_{\theta_{\star}} && }} 
 \end{equation}
 
 But by axioms {\bf LC1} and {\bf LC2} we have

 $$ \vcenter{\xymatrix@R=3ex@C=3ex{A_i \ar[rd]|{o}^u \ar[d]_{v'} \ar@/^5ex/[rrrrdd]|{o}^{\theta_{A_i}}_{\theta_{a} \Uparrow} \\
          C' \ar[rd]|{o}_{r} \ar@{}[r]|{\stackrel{\gamma}{\Rightarrow}} & C \ar[d]|<<<{o}^{s} \ar@/^2ex/[rrrd]|{o}^{\theta_{\star}} && \\
          & D  \ar[rrr]|{o}_{\theta_{\widetilde{\ast}}} & \ar@{}[ru]|{\theta_{d} \Uparrow} && E}}
          =
          \vcenter{\xymatrix@R=2.5ex@C=3ex{A_i \ar@/^5ex/[rrrrdd]|{o}^{\theta_{A_i}}_{\theta_{b'} \Uparrow} \ar[rd]_{v'} \\
          & C' \ar[rd]|{o}_{r}  \ar@/^2ex/[rrrd]|{o}^{\theta_{\star'}} & \ar@{}[rd]|{\theta_{c} \Uparrow} & \\
          && D  \ar[rr]|{o}_{\theta_{\widetilde{\ast}}} && E   }}$$
 
 %
%
 Composing this equation with $\theta_{a}^{-1}$ and $\theta_{c}^{-1} v'$ it follows \eqref{eq:toshow}, finishing the proof.
\end{proof}

\begin{remark} \label{rem:coroprooflemaposta}
 We note some facts that follow from the proof of Lemma \ref{lema:todofinitov2} and that will be used later. 
 We consider as in the proof of the lemma (Step 2) the set $S = \{\eta,\xi,...\}$ consisting of all the $\xi_i$ and all the compositions of the $\xi_i$ that appear in the given compositions, and consider for each $\eta$ in this family the premorphisms $\widetilde{\eta}$ as in \eqref{eq:premorfismostildeconeta}. Then:
  
  \begin{enumerate}
   \item If $\widetilde{\eta} = \widetilde{\xi}$, then $(\theta_a^{-1},\theta_{a'}^{-1},\theta_b,\theta_{b'}): \eta \Rightarrow \xi$ is an homotopy (this is immediate by the definition of homotopy).
   \item If $\eta \circ_{\gamma} \xi$ is in the set $S$, then $\widetilde{\eta \circ_{\gamma} \xi} = \widetilde{\eta} \circ \widetilde{\xi}$ (this is item B of the proof).
  \end{enumerate}
\end{remark}

\begin{remark}
 We will loosely say that we apply Lemma \ref{lema:todofinitov2} to a finite family of equations between compositions of premorphisms. We will omit to say that in this case we're considering all appearing premorphisms as the family $\xi_i$ and all the compositions appearing in the equations as the set $S = \{\eta,\xi,...\}$. 
\end{remark}

         \begin{\prop}[{cf. \cite[1.8]{DS}, vertical composition of homotopies}] \label{prop:vertcomphomot}
           The relation ``being homotopical'' is transitive, i.e. if $\xi_1 \sim \xi_2$ and $\xi_2 \sim \xi_3$ then $\xi_1 \sim \xi_3$.
         \end{\prop}
         
         \begin{proof}
         Apply Lemma \ref{lema:todofinitov2} to the equations $\xi_1 \sim \xi_2$ and $\xi_2 \sim \xi_3$. Then $\widetilde{\xi_1} = \widetilde{\xi_2} = \widetilde{\xi_3}$, which by remark \ref{rem:coroprooflemaposta}, item 1, implies that $\xi_1 \sim \xi_3$.
\end{proof}         

\begin{corollary} \label{coro:homotopicalesdeequiv}
The relation $\sim$ in Definition \ref{def:homotopical} is 
an equivalence relation between premorphisms that share the same domain and codomain. In fact, reflexivity and symmetry are obvious by definition, and transitivity follows from the previous proposition. \qed
\end{corollary}

\begin{\prop}[{cf. \cite[1.11]{DS}, horizontal composition of homotopies}] \label{prop:compbiendef}
 Consider composable premorphisms as follows $\xymatrix@C=8ex{ 
          {(x,\, A)} 
          \ar@<1ex>[r]^{\xi_1} 
          \ar@<-1ex>[r]_{\xi_2}
        & {(y,\, B)}
          \ar@<1ex>[r]^{\eta_1} 
          \ar@<-1ex>[r]_{\eta_2}
        & {(z,\, C)}
          }$. Then, given any two composites $\eta_1 \circ_{\gamma_1}  \xi_1$ and $\eta_2 \circ_{\gamma_2} \xi_2$, if $\xi_1 \sim \xi_2$ and $\eta_1 \sim \eta_2$ then 
          $\eta_1 \circ_{\gamma_1}  \xi_1  \sim \eta_2 \circ_{\gamma_2} \xi_2$.
\end{\prop}
\begin{proof}
We apply Lemma \ref{lema:todofinitov2} to the premorphisms $\{\xi_1,\xi_2,\eta_1,\eta_2\}$, the set \mbox{$S = \{\xi_1,\xi_2,\eta_1,\eta_2, \eta_1 \circ_{\gamma_1} \xi_1, \eta_2 \circ_{\gamma_2} \xi_2\}$} and the equations $\xi_1 \sim \xi_2$, $\eta_1 \sim \eta_2$. Then we have

$$\widetilde{\eta_1 \circ_{\gamma_1} \xi_1} = \widetilde{\eta_1} \circ \widetilde{\xi_1} = \widetilde{\eta_2} \circ \widetilde{\xi_2} = \widetilde{\eta_2 \circ_{\gamma_2} \xi_2},$$

\noindent where the equalities are justified by remark \ref{rem:coroprooflemaposta}, item 2 and by Lemma \ref{lema:todofinitov2}. By remark \ref{rem:coroprooflemaposta}, item 1, this implies that $\eta_1 \circ_{\gamma_1}  \xi_1  \sim \eta_2 \circ_{\gamma_2} \xi_2$.
\end{proof}

\begin{\de}[{cf. \cite[1.12,1.13]{DS}}] \label{olvidada1}
 Recall our Definition \ref{de:quasicatLF} of the quasicategory $\cc{L}_{\sigma}(F)$. 
 There is a category that we'll also denote by $\cc{L}_{\sigma}(F)$, whose arrows are the classes of premorphisms under the equivalence relation $\sim$. 
 We will denote the arrows of $\cc{L}_{\sigma}(F)$ with the same letter $\xi$, as it is clear from the context whether a premorphism or its class is being considered.
 
 
 Note that Proposition \ref{prop:compbiendef} shows that the composition of classes is well defined and may be computed over an arbitrary $2$-cell $\gamma$ as in \eqref{eq:compdepremorf}. Identities are $(x,A) \xr{(id_A,id_x,id_A)} (x,A)$ and composition is associative since we can choose the $2$-cells $\gamma_i$, $i=1,2,3$ as in \eqref{eq:paraasociat}.
\end{\de}

\begin{\prop}[{cf. \cite[Proof of 2.4 (b)]{DS}}] \label{prop:mismoalphaybeta}
 Consider a finite sequence of functors $\cc{A} \mr{F_k} \Cat$ and premorphisms of $\cc{L}_{\sigma}(F_k)$ as follows  
$\xymatrix@ur@R=3.5ex@C=3.5ex
         {
        F_k \ar[r]^{x_k}  \ar@{}[dr]|{\xi_k \Downarrow}  \ar[d]_{y_k} 
          & A \ar[d]|{o}^{u}
          \\
            B  \ar[r]_{v} 
          & C 
         }
, \;
\xymatrix@ur@R=3.5ex@C=3.5ex
         {
        F_k \ar[r]^{x_k}  \ar@{}[dr]|{\eta_k \Downarrow}  \ar[d]_{y_k} 
          & A \ar[d]|{o}^{s}
          \\
            B  \ar[r]_{t} 
          & D 
         }.$
         If $\eta_k \sim \xi_k$ for each $k$, then there are 2-cells $\alpha_1$, $\alpha_2$, $\beta_1$, $\beta_2$ which define all the
homotopies  $(\alpha_1, \alpha_2, \beta_1, \beta_2): \eta_k \Rightarrow \xi_k$.
\end{\prop}

\begin{proof}
We apply Lemma \ref{lema:todofinitov2} to the equations $\eta_k \sim \xi_k$. In this way we have $\widetilde{\eta_k} = \widetilde{\xi_k}$, but since all the $\xi_k$ share the same $(u,v)$ then all the $\widetilde{\xi_k}$ are constructed from the $\xi_k$ by pasting the same $2$-cells $\mu,\nu$ that we denote by $\alpha_2, \beta_2$ respectively. Similarly, all the $\widetilde{\eta_k}$ are constructed from the $\eta_k$ by pasting the same $2$-cells that we denote $\alpha_1,\beta_1$. Then by remark \ref{rem:coroprooflemaposta}, item 1, $(\alpha_1,\alpha_2,\beta_1,\beta_2) : \eta_k \Rightarrow \xi_k$ for each $k$ as desired.
\end{proof}

\begin{theorem}[{cf. \cite[1.19]{DS}}] \label{olvidada2}
 The following formulas define a $\sigma$-cone $F \Mr{\lambda} \cc{L}_{\sigma}(F)$ which is the $\sigma$-colimit of $F$ in $\Cat$. For $A \in \A$, $x \in FA$, $x \mr{\xi} y$ in $FA$, $A \mr{u} B$ in $\A$
 $$
\lambda_A(x) = F \mr{x} A \,, 
\hspace{3ex}
\lambda_A(\xi) =
\xymatrix@ur@R=4ex@C=4ex
         {
            F \ar[r]^x  \ar@{}[dr]|{\xi \Downarrow}  \ar[d]_y 
          & A \ar[d]^{{id}}
          \\
            A  \ar[r]|{o}_{{id}} 
          & A
         },
\hspace{3ex}
(\lambda_u)_x = 
         \xymatrix@ur@R=4ex@C=4ex
         {
            F \ar[r]^{ux}  \ar@{}[dr]|{id \Downarrow}  \ar[d]_x 
          & B \ar[d]|{o}^{{id}}
          \\
            A  \ar[r]_{{u}} 
          & B
         } $$
%
\end{theorem}
\begin{proof}
 %
 It is immediate to show that $\lambda$ is a functor, that $\lambda_u$ are natural transformations and axioms {\bf LC0}, {\bf LC1} in Definition \ref{def:sigmacone}. 
 If $u \in \Sigma$, then for each $x$ we have the premorphisms $(\lambda_u^{-1})_x = 
         \xymatrix@ur@R=4ex@C=4ex
         {
            F \ar[d]_{ux}  \ar@{}[dr]|{id \Downarrow}  \ar[r]^x 
          & B \ar[d]|{o}^{{u}} 
          \\
            A \ar[r]_{{id}} 
          & B
         },$
%
         showing that $\lambda_u$ is invertible.
         
To prove axiom {\bf LC2}, given $A \cellrd{u}{\gamma}{v} B$, we have to show that

$$\vcenter{\xymatrix@R=1ex@C=4ex
         { & B \ar[rd]^{id} \\
          & {\txt{\scriptsize{$(F\gamma)_x \Downarrow$}}} & B \ar[rd]^{id} \\
           F \ar[r]^{vx} \ar[rdd]_{x} \ar[ruu]^{ux} & B \ar[ru]|{id} \ar[rd]|{id} & {\txt{\scriptsize{$id \Downarrow$}}} & B \\
          & {\txt{\scriptsize{$id \Downarrow$}}} & B \ar[ru]_{id} \\
          & A \ar[ru]_{v} }}  
         \sim
         \vcenter{\xymatrix@R=1ex@C=4ex
         { & B \ar[rd]^{id} \\
         F \ar[rd]_{x} \ar[ru]^{ux} & {\txt{\scriptsize{$id \Downarrow$}}} & B \\
         & A \ar[ru]_{u} }}          $$

 
 We consider a $\sigma$-cone $\theta$ with arrows in $\Sigma$ with vertex $C$ for the diagram determined by the $2$-cell $\gamma$ (recall Proposition \ref{prop:filtsiicone}). Then a simple computation shows that
 
          $$ 
         \vcenter{\xymatrix@R=0.2ex@C=3ex
         {
          & B \ar[rdd]|{o}^{id} \ar@/^2ex/[rrrdd]|{o}^{\theta_B}
          \\
          & {\txt{\scriptsize{$(F\gamma)_x \!\Downarrow$}}} & {\;\; \txt{\scriptsize{$id \!\Downarrow$}}} 
          \\
          F \ar[ruu]^{ux} \ar[rdd]_{x} \ar[rr]|{vx}
          && B \ar[rr]|{o}^{\theta_B}
          && C
          \\
          &{\txt{\scriptsize{$id \!\Downarrow$}}} & {\txt{\;\; \scriptsize{$\theta_v \!\Downarrow$}}}
          \\
          & A \ar[ruu]_{v} \ar@/_2ex/[rrruu]|{o}_{\theta_A}
         }}
         =
   \vcenter{\xymatrix@R=0.2ex@C=3ex
         {
          & B \ar[rdd]|{o}^{id} \ar@/^2ex/[rrrdd]|{o}^{\theta_B}
          \\
          && {\;\; \txt{\scriptsize{$id \!\Downarrow$}}} 
          \\
           F \ar[ruu]^{ux} \ar[rdd]_{x}
          & {\txt{\scriptsize{$id \!\Downarrow$}}}
          & B \ar[rr]|{o}^{\theta_B}
          && C
          \\
          && {\txt{\;\; \scriptsize{$\theta_u \!\Downarrow$}}}
          \\
          & A \ar[ruu]_{u} \ar@/_2ex/[rrruu]|{o}_{\theta_A}
         }}$$
 
 We will now show the universal property. Note that since $\sigma$-colimits exist in $\Cat$ \mbox{\cite[I,7.9.1 iii)]{GRAY},} \cite{DDS}, 
 it suffices to show the $1$-dimensional property (see \cite[(3.4) and below]{K1}). Then if we consider $F \Mr{h} \cc{X}$ a $\sigma$-cone, we have to prove that there is a unique functor $\cc{L}_{\sigma}(F) \mr{\widetilde{h}} \cc{X}$ such that $\widetilde{h} \lambda = h$. The only possible definition on objects is $\widetilde{h}(x,A) = h_A(x)$, and given any premorphism $\xymatrix@R=4ex@C=4ex@ur
         {
            F \ar[r]^x  \ar@{}[dr]|{\xi \Downarrow}  \ar[d]_y 
          & A \ar[d]|{o}^u
          \\
            B  \ar[r]_v 
          & C 
         }$, the equality $\xi = \lambda_v \circ \xi \circ \lambda_u^{-1}$ forces the definition $\widetilde{h}(\xi) = \vcenter{\xymatrix@R=0.2ex@C=3ex
         {
          & A \ar[rdd]|{o}_u \ar@/^2ex/[rrrdd]^{h_A}
          \\
          && {\;\; \txt{\scriptsize{$h_{u}^{-1} \!\Downarrow$}}} 
          \\
          F \ar[ruu]^x \ar[rdd]_y
          & {\txt{\scriptsize{$\xi \!\Downarrow$}}}
          & C \ar[rr]^{h_C}
          && \cc{X}
          \\
          && {\txt{\;\; \scriptsize{$h_v \!\Downarrow$}}}
          \\
          & B \ar[ruu]^v \ar@/_2ex/[rrruu]_{h_B} 
         }}$       

 To show that this is well defined, consider a homotopy $(\alpha_1,\alpha_2,\beta_1,\beta_2): \xi_1 \Rightarrow \xi_2$, 
 then we compute
 
 $$\vcenter{\xymatrix@R=0.2ex@C=3ex{
          & A \ar[rdd]|{o}_{u_1} \ar@/^2ex/[rrrdd]^{h_A}
          \\
          && {\;\; \txt{\scriptsize{$h_{u_1}^{-1} \!\Downarrow$}}} 
          \\
          F \ar[ruu]^x \ar[rdd]_y
          & {\txt{\scriptsize{$\xi_1 \!\Downarrow$}}}
          & C_1 \ar[rr]^{h_{C_1}}
          && \cc{X}           \\
          && {\txt{\;\; \scriptsize{$h_{v_1} \!\Downarrow$}}}           \\
          & B \ar[ruu]^{v_1} \ar@/_2ex/[rrruu]_{h_B} }}
         \stackrel{(*)}{=}
   \vcenter{\xymatrix@R=0.2ex@C=3ex
         {
          & A \ar[rdd]|{o}_{u_1} \ar@/^2ex/[rrrdd]|>>>>{o}^>>>>{w_A} \ar@/^5ex/[rrrrdd]^{h_A}_{h_{w_A}^{-1} \Downarrow}
          \\
          && {\;\; \txt{\scriptsize{$\alpha_1 \!\Downarrow$}}} 
          \\
          F \ar[ruu]^{x} \ar[rdd]_{y}
          & {\txt{\scriptsize{$\xi_1 \!\Downarrow$}}}
          & C_1 \ar[rr]|{o}^{w_1}
          && D \ar[r]^{h_D} & \cc{X}
          \\
          && {\txt{\;\; \scriptsize{$\beta_1 \!\Downarrow$}}}
          \\
          & B \ar[ruu]^{v_1} \ar@/_2ex/[rrruu]|>>>>{o}_>>>>{w_B} \ar@/_5ex/[rrrruu]_{h_B}^{h_{w_B} \Downarrow}
         }}=$$
         
         $$\vcenter{\xymatrix@R=0.2ex@C=3ex
         {
          & A \ar[rdd]|{o}_{u_2} \ar@/^2ex/[rrrdd]|>>>>{o}^>>>>{w_A} \ar@/^5ex/[rrrrdd]^{h_A}_{h_{w_A}^{-1} \Downarrow}
          \\
          && {\;\; \txt{\scriptsize{$\alpha_2 \!\Downarrow$}}} 
          \\
          F \ar[ruu]^{x} \ar[rdd]_{y}
          & {\txt{\scriptsize{$\xi_2 \!\Downarrow$}}}
          & C_2 \ar[rr]|{o}^{w_2}
          && D \ar[r]^{h_D} & \cc{X}
          \\
          && {\txt{\;\; \scriptsize{$\beta_2 \!\Downarrow$}}}
          \\
          & B \ar[ruu]^{v_2} \ar@/_2ex/[rrruu]|>>>>{o}_>>>>{w_B} \ar@/_5ex/[rrrruu]_{h_B}^{h_{w_B} \Downarrow}
         }}
         \stackrel{(*)}{=}
         \vcenter{\xymatrix@R=0.2ex@C=3ex
         {
          & A \ar[rdd]|{o}_{u_2} \ar@/^2ex/[rrrdd]^{h_A}
          \\
          && {\;\; \txt{\scriptsize{$h_{u_2}^{-1} \!\Downarrow$}}} 
          \\
          F \ar[ruu]^x \ar[rdd]_y
          & {\txt{\scriptsize{$\xi_2 \!\Downarrow$}}}
          & C_2 \ar[rr]^{h_{C_2}}
          && \cc{X}
          \\
          && {\txt{\;\; \scriptsize{$h_{v_2} \!\Downarrow$}}}
          \\
          & B \ar[ruu]^{v_2} \ar@/_2ex/[rrruu]_{h_B} 
         }} $$
         
         The equalities marked with an $(*)$ are justified from the following equality and the corresponding equalities for $\beta_1$, $\alpha_2$ and $\beta_2$, all of which follow immediately from axioms {\bf LC1} and {\bf LC2}.
         
          $$ \vcenter{\xymatrix@R=1.2ex@C=3ex
         {
           A \ar[rdd]_{u_1} \ar@/^4ex/[rrrdd]^{h_A} \ar@{}[rrrd]_{h_{u_1}^{-1} \Downarrow}
          \\
          &&& 
          \\
          & C_1 \ar[rr]_{h_{C_1}} && \cc{X}  }}
          =
         \vcenter{\xymatrix@R=0.2ex@C=3ex
         {
           A \ar[rdd]_{u_1} \ar@/^2ex/[rrrdd]|>>>>{o}^>>>>{w_A} \ar@/^5ex/[rrrrdd]^{h_A}_{h_{w_A}^{-1} \Downarrow}
          \\
          & {\;\; \txt{\scriptsize{$\alpha_1 \!\Downarrow$}}} 
          \\
          & C_1 \ar[rr]|{o}^{w_1} \ar@/_4ex/[rrr]_{h_{C_1}}^{h_{w_1} \Downarrow}
          && D \ar[r]^{h_D} & \cc{X}  }}
$$

 \end{proof}

\section{A basic exactness property of \boldmath $\cc{C}at$} \label{sub:exactness}
 The following theorem, 
 whose proof will occupy this whole section,
 generalizes \cite[Theorem 2.4]{DS} not only to the $\sigma$-case but also to all finite weighted bilimits instead of cotensors. The case where $\Sigma$ consists of all the arrows of $\A$ was considered in \cite{Canevali}.

 
 We refer to \cite{K2}, \cite{S2}, \cite{DDS} among other choices for the definition of the bilimit of a $2$-functor $\A \mr{F} \B$ weighted by a $2$-functor $\A \mr{W} \Cat$.
%
%
%
 We consider below a simple version of \emph{finite weight} which is sufficient for our purposes.
 \begin{\de} \label{def:finitelimit} 
 We say that a $2$-functor $\A \mr{W} \Cat$ is a \emph{finite weight} if $\A$ is a finite $2$-category and for each $A \in \A$, $WA$ is a finite category.  Given $2$-functors $\A \mr{W} \Cat$, $\A \mr{F} \B$ the bilimit ${}_{bi} \! \{W,F\}$ is said to be a \emph{finite bilimit} if $W$ is a finite weight.
 \end{\de}
 
\begin{theorem} \label{teo:conmutan}
$\sigma$-filtered $\sigma$-colimits commute with finite (weighted) bilimits in $\Cat$. This means precisely that, for a finite $2$-category $\cc{I}$, and $2$-functors $\A \times \cc{I} \mr{F} \Cat$, $\cc{I} \mr{W} \Cat$, if $W$ is a finite weight, then
the canonical comparison functor 

\begin{equation} \label{eq:comparisonfunctor}
\Diamond: \coLim{A \in \A}{ \;{}_{bi} \!\{ W,F(A,-) \} } \mr{\approx} \;{}_{bi} \!\{ W,\coLim{A \in \A}{F(A,-)} \} 
\end{equation}

\noindent is an equivalence of categories.
Observe that the bilimits can be computed in $\Cat$ as pseudo-limits. Despite that, commutativity holds not up to isomorphism but only up to equivalence, i.e. in the ``bi-sense''.  
\end{theorem}
 As it is known, weighted bilimits can be constructed using bicotensors, biproducts and biequalizers (see \cite{SCorr}). These bilimits exist in $\Cat$ in the stronger pseudo form. Furthermore, the proof in \cite{SCorr} also yields the \emph{finite} case, see \cite[Corollary 6.12]{Canevali} for a detailed proof of the following: \emph{finite weighted bilimits in $\Cat$ can be constructed (as pseudo-limits) from finite cotensors, finite $2$-products and pseudo-equalizers}. Therefore to prove Theorem \ref{teo:conmutan} it suffices to check commutativity for these three special kind of limits, which we do in the next three propositions. 
  
 
\vspace{1ex}

The reader may notice that the proofs of these three propositions are somewhat similar, and in fact a direct unified proof of Theorem \ref{teo:conmutan} is also possible. This proof is obtained by describing the categories appearing in \eqref{eq:comparisonfunctor} and the action of $\Diamond$, using the constructions of $\sigma$-filtered $\sigma$-colimits in $\Cat$, and the construction of weighted bilimits in $\Cat$. Since the weight $W$ takes values in finite categories, Lemma \ref{lema:todofinitov2} can also be used to prove this general case, but since some computations become much harder than the ones in the particular cases of the three propositions below we have refrained from doing so.

\vspace{1ex}

 When using Lemma \ref{lema:todofinitov2}, we will denote the arrows $w_i, z_i, t_i$ in \eqref{eq:premorfismostilde} by $A \mr{w_A} E$, $B \mr{w_B} E$, and so on, this is justified by item 1 in Lemma \ref{lema:todofinitov2}.

\vspace{1ex}
 
 \noindent {\bf Finite cotensors.} Note that if $C, D \in \Cat$, we have $\{D, C\}=C^D$ the category of functors from $C$ to $D$. 
 Given $F: \A \mr{} \Cat$, and $W$ a category, since cotensors are computed pointwise we have the $2$-functor $\{W,F\}(A) = (FA)^{W}$, we will denote $F^W = \{W,F\}$. We have a comparison functor $\Diamond: \cc{L}_{\sigma}(F^W) \mr{} \cc{L}_{\sigma}(F)^W$, we describe its action on objects and on arrows. 
 
 \begin{enumerate}
 \item An object in $\cc{L}_{\sigma}(F^W)$ is a pair $(x,A)$, where $x \in F^W(A) = (FA)^{W}$. This is given by $F \mr{x_k} A$, $k \in W$, and arrows $x_p \mr{x_f} x_q$, $p \mr{f} q \in W$, satisfying equations $x_{id_k}=id_{x_k}$ for all $k \in W$, $ x_{f \circ g} = x_f \circ x_g$ 
for all composable pairs $f,\; g \in W$.
 
 \item An object in $\cc{L}_{\sigma}(F)^W$ is a functor $W \xr{(x,\varphi)} \cc{L}_{\sigma}(F)$ given explicitly by:
$$
F \mr{x_k} A_k ,\; k \in W, \hspace{2ex} 
\xymatrix@R=4ex@C=4ex@ur
         {
            F \ar[r]^{x_p}  \ar@{}[dr]|{\varphi_f \Downarrow}
	    \ar[d]_{x_q}  
          & A_p \ar[d]|{o}^{u_f}
          \\
            A_q  \ar[r]_{v_f} 
          & A_f 
         },
\hspace{2ex} p \mr{f} q  \in  W,
$$
\noindent satisfying equations $\varphi_{id_k}\,\sim\, id_{x_k}$ for all $k \in W$, $ \varphi_{f \circ g} \,\sim\, \varphi_f \circ \varphi_g$ 
for all composable pairs $f,\; g$. 

  \item To describe $\Diamond(x,A)$, consider the object in $\cc{L}_{\sigma}(F)^W$ with $F \mr{x_k} A$, $\varphi_f = \lambda_A(x_f)$.
%
%

 \item A premorphism $\vcenter{\xymatrix@R=4ex@C=4ex@ur
         {
            F^{W} \ar[r]^{x}  \ar@{}[dr]|{\xi \Downarrow}
	    \ar[d]_{y}  
          & A \ar[d]|{o}^{u}
          \\
            B  \ar[r]_{v} 
          & C 
         }}$ in $\cc{L}_{\sigma}(F^W)$ is a natural transformation 
         $\vcenter{\xymatrix@R=2ex@C=2ex{ & FA \ar[rd]^{F(u)} \\
         W \ar[rd]_y \ar[ru]^x & {\txt{\scriptsize{$\xi \!\Downarrow$}}} & FC \\
         & FB \ar[ru]_{F(v)} }}$. The naturality equations are, for $p \mr{f} q$ in $W$, $Fv(y_f) \circ \xi_p = \xi_q \circ Fu(x_f)$.
         
Two premorphisms $\xi_1$, $\xi_2$ are equivalent if there are $2$-cells $(\alpha_1,\alpha_2,\beta_1,\beta_2)$, with $\alpha_i$ invertible $i=1,2$, such that
         
         $$ 
         \vcenter{\xymatrix@R=0.2ex@C=3ex
         {
          & A \ar[rdd]|{o}_{u_1} \ar@/^2ex/[rrrdd]|{o}^{w_A}
          \\
          && {\;\; \txt{\scriptsize{$\alpha_1 \!\Downarrow$}}} 
          \\
          F^{W} \ar[ruu]^{x} \ar[rdd]_{y}
          & {\txt{\scriptsize{$\xi_1 \!\Downarrow$}}}
          & C_1 \ar[rr]|{o}^{w_1}
          && D
          \\
          && {\txt{\;\; \scriptsize{$\beta_1 \!\Downarrow$}}}
          \\
          & B \ar[ruu]^{v_1} \ar@/_2ex/[rrruu]|{o}_{w_B}
         }}
         =
         \vcenter{\xymatrix@R=0.2ex@C=3ex
         {
          & A \ar[rdd]|{o}_{u_2} \ar@/^2ex/[rrrdd]|{o}^{w_A}
          \\
          && {\;\; \txt{\scriptsize{$\alpha_2 \!\Downarrow$}}} 
          \\
          F^{W} \ar[ruu]^{x} \ar[rdd]_{y}
          & {\txt{\scriptsize{$\xi_2 \!\Downarrow$}}}
          & C_2 \ar[rr]|{o}^{w_2}
          && D
          \\
          && {\txt{\;\; \scriptsize{$\beta_2 \!\Downarrow$}}}
          \\
          & B \ar[ruu]^{v_2} \ar@/_2ex/[rrruu]|{o}_{w_B}
         }}$$
         
         \item An arrow $(x,\varphi) \mr{\xi} (y,\psi)$ of $\cc{L}_{\sigma}(F)^W$ is a natural transformation between the functors, it consists of a family of premorphisms 
         $\vcenter{\xymatrix@R=4ex@C=4ex@ur
         {
            F \ar[r]^{x_k}  \ar@{}[dr]|{\xi_k \Downarrow}
	    \ar[d]_{y_k}  
          & A_k \ar[d]|{o}^{u_k}
          \\
            B_k  \ar[r]_{v_k} 
          & C_k 
         }}$ satisfying the naturality equations, for $p \mr{f} q$ in $W$, $\psi_f \circ \xi_p \sim \xi_q \circ \varphi_f$.
         
         \item The description of $\Diamond$ on arrows is $(\Diamond(\xi))_k = \xi_k$ for all $k \in W$.
\end{enumerate}
 
\begin{proposition}[{cf. \cite[2.4]{DS}}] \label{teo:conmutanconcotensors}
 $\sigma$-filtered $\sigma$-colimits commute with finite cotensors in $\Cat$. This means precisely that, for a finite category $W$, the canonical functor 
 \mbox{$\Diamond: \cc{L}_{\sigma}(F^W) \mr{} \cc{L}_{\sigma}(F)^W$} is an equivalence of categories.
\end{proposition}
\begin{proof}
 We follow the structure of the proof in \cite{DS}, but by virtue of previous lemmas the proof here (particularly of item (a)) is simpler. 
 
 \vspace{1ex}
 
 \noindent \textit{(a) essentially surjective.} 
Given an object $(x, \varphi)$ in $\cc{L}_{\sigma}(F)^W$ as above, we apply Lemma \ref{lema:todofinitov2} to have an object $A$, arrows $A_k \mr{w_k} A$, and premorphisms $\widetilde{\varphi_f} \sim \varphi_f$, $\xymatrix@R=4ex@C=4ex@ur
         {
            F \ar[r]^{x_p}  \ar@{}[dr]|{\widetilde{\varphi_f} \Downarrow}
	    \ar[d]_{x_q}  
          & A_p \ar[d]|{o}^{w_p}
          \\
            A_q  \ar[r]|{o}_{w_q} 
          & A 
         }$, such that
the equations $\widetilde{\varphi_{f \circ g}} \,=\, \widetilde{\varphi_f} \circ \widetilde{\varphi_g}$, $\widetilde{\varphi_{id_k}} = id_{w_k x_k}$ hold in $FA$. These data amount to a 
factorization (up to isomorphism) 
$$
\xymatrix
        {
         & W \ar[d]^{(x,\varphi)} \ar[ld]_{y}^{\;\;\cong}
         \\
         FA \ar[r]^{\lambda_A}
         & \cc{L}_{\sigma}(F)
        } 
$$
\noindent where the functor $W \mr{y} FA$ is defined by the formulas $y_k = w_k x_k$, $y_f = \widetilde{\varphi_f}$, and
the identity $2$-cells yield the isomorphism $(x,\varphi) \cong \Diamond (y,A)$. 

%
%
%
%

\vspace{1ex}

 \noindent \textit{(b) faithful.} 
 If $\xi$, $\eta \in \cc{L}_{\sigma}(F^W)$ are such that $\Diamond (\xi)=\Diamond (\eta)$, then $\xi_k=\eta_k$ for all $k\in W$. From Proposition \ref{prop:mismoalphaybeta}, this implies that $\xi\,\sim\, \eta$ which concludes the proof.

\vspace{1ex}

 \noindent \textit{(c) full:} An arrow $(x,\lambda_A(x_f)) \mr{\xi} (y,\lambda_B(y_f))$ is a family of premorphisms 
         $\vcenter{\xymatrix@R=4ex@C=4ex@ur
         {
            F \ar[r]^{x_k}  \ar@{}[dr]|{\xi_k \Downarrow}
	    \ar[d]_{y_k}  
          & A \ar[d]|{o}^{u_k}
          \\
            B  \ar[r]_{v_k} 
          & C_k 
         }}$ satisfying the naturality equations, for $p \mr{f} q$ in $W$, $\lambda_B(y_f) \circ \xi_p \sim \xi_q \circ \lambda_A(x_f)$. We apply Lemma \ref{lema:todofinitov2} to have $C$ and equivalent premorphisms such that the equations 
\mbox{$\widetilde{\lambda_B(y_f)} \circ \widetilde{\xi_p} = \widetilde{\xi_q} \circ \widetilde{\lambda_A(x_f)}$} hold. We claim that $\vcenter{\xymatrix@R=4ex@C=4ex@ur
         {
            F^{W} \ar[r]^{x}  \ar@{}[dr]|{\widetilde{\xi} \Downarrow}
	    \ar[d]_{y}  
          & A \ar[d]|{o}^{w_A}
          \\
            B  \ar[r]|{o}_{w_B} 
          & C
         }}$ is a natural transformation. This is so because by 
         item 2.i) in Lemma \ref{lema:todofinitov2},
         we have $\widetilde{\lambda_A(x_f)} = F(w_A)(x_f)$, $\widetilde{\lambda_B(y_f)} = F(w_B)(y_f)$, thus the naturality equations hold.
%
%
\end{proof}
%
%
 \noindent {\bf $2$-products.} 
 The empty product offers no difficulty and we leave it to the reader. Given $F,G: \A \mr{} \Cat$, we have a comparison functor $\cc{L}_{\sigma}(F \times G) \mr{\Diamond} \cc{L}_{\sigma}(F) \times \cc{L}_{\sigma}(G)$, given by the formulas $\Diamond(F \times G \xr{(x,y)} A) = (F \mr{x} A, G \mr{y} A)$, $\Diamond(\xi,\eta) = (\xi,\eta)$. The reader can look at \cite[\S 7.3]{Canevali} for 
 a detailed description of $\cc{L}(F \times G) \mr{\Diamond} \cc{L}(F) \times \cc{L}(G)$, easily modifiable to our case where a subcategory $\Sigma$ of arrows is present, though our proof of the following theorem differs substantially from the proof of \cite[Theorem 7.10]{Canevali}, as we have available the powerful lemmas from section \ref{sub:DSC}.
 %
%
\begin{proposition} \label{teo:conmutaconproductos}
 $\sigma$-filtered $\sigma$-colimits commute with finite $2$-products in $\Cat$. This means precisely that the canonical functor $\cc{L}_{\sigma}(F \times G) \mr{\Diamond} \cc{L}_{\sigma}(F) \times \cc{L}_{\sigma}(G)$ is an equivalence of categories.
\end{proposition}
 
\begin{proof}
 \noindent \textit{(a) essentially surjective.} Given an object $(F \mr{x} A,G \mr{y} B)$ of $\cc{L}_{\sigma}(F) \times \cc{L}_{\sigma}(G)$, by axiom $\sigma F0$ 
 we have $\xymatrix{A \ar[r]|{o}^{w_A} & E}$, $\xymatrix{B \ar[r]|{o}^{w_B} & E}$. Then 
 $$(F \mr{x} A,G \mr{y} B) \cong (F \xr{F(w_A)(x)} E,G \xr{G(w_B)(y)} E) = \Diamond(F \times G \xr{(F(w_A)(x),G(w_B)(y))} E)$$
%
%
\noindent \textit{(b) faithful.}  It is immediate from Proposition \ref{prop:mismoalphaybeta}.

\vspace{1ex}

 \noindent \textit{(c) full.} Consider a premorphism of $\cc{L}_{\sigma}(F) \times \cc{L}_{\sigma}(G)$ between two objects in the image of $\Diamond$, that is 

$$
( \xymatrix@ur@-.5pc{F \ar[r]^x \ar[d]_{x'} \ar@{}[dr] |{\Downarrow \xi} & A \ar[d]|{o}^{u} \\ A' \ar[r]_{v} & C},
\xymatrix@ur@-.5pc{G \ar[r]^y \ar[d]_{y'} \ar@{}[dr] |{\Downarrow \eta} & A \ar[d]|{o}^r \\ A' \ar[r]_{s} & D} 
): (F \mr{x} A,G \mr{y} A) \mr{} (F \mr{x'} A',G \mr{y'} A').$$
 
 We apply Lemma \ref{lema:todofinitov2} to have equivalent premorphisms 
 $\xymatrix@ur@-.5pc{F \ar[r]^x \ar[d]_{x'} \ar@{}[dr] |{\Downarrow \widetilde{\xi}} & A \ar[d]|{o}^{w_A} \\ A' \ar[r]|{o}_{w_{A'}} & E},
\xymatrix@ur@-.5pc{G \ar[r]^y \ar[d]_{y'} \ar@{}[dr] |{\Downarrow \widetilde{\eta}} & A \ar[d]|{o}^{w_A} \\ A' \ar[r]|{o}_{w_{A'}} & E}$. 
Then $(\xi,\eta) \sim (\widetilde{\xi},\widetilde{\eta}) = \Diamond(\widetilde{\xi},\widetilde{\eta})$.
\end{proof}


 \noindent {\bf Pseudo-equalizers.} 
 We will show that $\sigma$-filtered $\sigma$-colimits commute with biequalizers in $\Cat$, this is the $\sigma$-version of \cite[Theorem 7.19]{Canevali}, we give a shorter complete proof. 
 Given 2-functors and 2-natural transformations $\xymatrix{\mathcal{A} \ar@/^/@<1ex>[r]^G \ar@/_/@<-0.5ex>[r]_H \ar@{} [r] |{{\Downarrow \alpha\: \Downarrow \beta}} & **[r]\Cat}$, we consider its pseudo-equalizer in $\cc{H}om_s(\A,\Cat)$, which is constructed pointwise, then we have a $2$-functor $P: \A \mr{} \Cat$. 
  We also consider the induced diagram $\xymatrix{\cc{L}_{\sigma}(G) \ar@<1ex>[r]^{\alpha} \ar@<-0.5ex>[r]_\beta & \cc{L}_{\sigma}(H)}$ and construct its pseudo-equalizer $Eq$ in $\Cat$, then we have a comparison functor $\cc{L}_{\sigma}(P) \mr{\Diamond} Eq$. 
As with $2$-products, we refer the reader to \cite[\S 7.4]{Canevali} for more details on the description of $\Diamond$ below:
 
 \begin{enumerate}
  \item  Objects of $\cc{L}_{\sigma}(P)$ are $P \xr{(x,y,\gamma,\delta)} A$ with $G \mr{x} A$, $H \mr{y} A$, and $\gamma, \delta$ isomorphisms in $HA$, $\alpha_A(x) \mr{\gamma} y$, $\beta_A(x) \mr{\delta} y$. 
  \item  Objects of $Eq$ are $(x,y,\varphi, \psi)$, where $G \mr{x} A$, $H \mr{y} B$, but now $\varphi, \psi$ are isomorphisms in $\cc{L}_{\sigma}(H)$, $(\alpha_A(x),A) \mr{\varphi} (y,B)$, $(\beta_A(x),A) \mr{\psi} (y,B)$.
  \item The description of $\Diamond$ on objects is $\Diamond(x,y,\gamma,\delta) = (x,y,\lambda_A(\gamma),\lambda_A(\delta))$, where \mbox{$HA \mr{\lambda_A} \cc{L}_{\sigma}(H)$} is the inclusion in the $\sigma$-colimit. Note that, since \mbox{$\lambda_A(\gamma) = (id,\gamma,id)$} and we usually denote a premorphism $(u,\xi,v)$ just by $\xi$, we may write \mbox{$\Diamond(x,y,\gamma,\delta) = (x,y,\gamma,\delta)$,} but we must keep in mind the distinction between premorphisms of $\cc{L}_{\sigma}(H)$ and arrows in $HA$.
  \item A premorphism of $\cc{L}_{\sigma}(P)$ between objects $P \xr{(x,y,\gamma,\delta)} A$, $P \xr{(x',y',\gamma',\delta')} A'$ is given by a 
 pair of premorphisms
 $\xymatrix@ur@-.5pc{G \ar[r]^x \ar[d]_{x'} \ar@{}[dr] |{\Downarrow \xi} & A \ar[d]|{o}^u \\ A' \ar[r]_{v} & C},
\xymatrix@ur@-.5pc{H \ar[r]^y \ar[d]_{y'} \ar@{}[dr] |{\Downarrow \eta} & A \ar[d]|{o}^u \\ A' \ar[r]_{v} & C}$ in $\cc{L}_{\sigma}(G)$ and $\cc{L}_{\sigma}(H)$ 
 satisfying the following two equations in $HC$
 \begin{align*} \eta \circ H(u)(\gamma) &= H(v)(\gamma') \circ \alpha_C(\xi)\\
 \eta \circ H(u)(\delta) &= H(v)(\delta') \circ \beta_C(\xi)
\end{align*}

Two premorphisms $(\xi_1,\eta_1)$ and $(\xi_2,\eta_2)$ are equivalent if and only if there are homotopies given by the same $2$-cells $(\alpha_1,\alpha_2,\beta_1,\beta_2) : \xi_1 \Rightarrow \xi_2$ and $(\alpha_1,\alpha_2,\beta_1,\beta_2) : \eta_1 \Rightarrow \eta_2$.

\item A premorphism of $Eq$ between $(x,y,\varphi, \psi)$, $(x',y',\varphi', \psi')$ is also given by a 
 pair of premorphisms
 $\xymatrix@ur@-.5pc{G \ar[r]^x \ar[d]_{x'} \ar@{}[dr] |{\Downarrow \xi} & A \ar[d]|{o}^{u} \\ A' \ar[r]_{v} & C},
\xymatrix@ur@-.5pc{H \ar[r]^y \ar[d]_{y'} \ar@{}[dr] |{\Downarrow \eta} & B \ar[d]|{o}^r \\ B' \ar[r]_{s} & D}$ in $\cc{L}_{\sigma}(G)$ and $\cc{L}_{\sigma}(H)$, 
 but now the following similar two equations are satisfied up to homotopy, i.e. in $\cc{L}_{\sigma}(H)$
 \begin{align*} \eta \circ \varphi \sim \varphi' \circ \alpha_C(\xi)\\
 \eta \circ \psi \sim \psi' \circ \beta_C(\xi)
\end{align*}
Two premorphisms $(\xi_1,\eta_1)$ and $(\xi_2,\eta_2)$ are equivalent if and only if $\xi_1 \sim \xi_2$ in $\cc{L}_{\sigma}(G)$ and $\eta_1 \sim \eta_2$ in $\cc{L}_{\sigma}(H)$.
  \item The description of $\Diamond$ on arrows is $\Diamond(\xi,\eta) = (\xi,\eta)$. 
 \end{enumerate}

\begin{proposition} \label{teo:conmutaconbieq}
$\sigma$-filtered $\sigma$-colimits commute with biequalizers in $\Cat$. This means precisely that the canonical functor $\cc{L}_{\sigma}(P) \mr{\Diamond} Eq$ is an equivalence of categories.
\end{proposition}

\begin{proof}
\noindent \textit{(a) essentially surjective.} Given an object $(x,y,\varphi,\psi)$ of $Eq$ as above, by Lemma \ref{lema:todofinitov2} (with no equations involved) we have
$\xymatrix{A \ar[r]|{o}^{w_A} & E}$, $\xymatrix{B \ar[r]|{o}^{w_B} & E}$ and premorphisms $
\xymatrix@-.5pc@ur
         {
            H \ar[r]^{\alpha_A(x)}  \ar@{}[dr]|{\widetilde{\varphi} \Downarrow}  \ar[d]_y 
          & A \ar[d]|{o}^{w_A}
          \\
            B  \ar[r]|{o}_{w_B}
          & E 
         }$;
$\xymatrix@-.5pc@ur
         {
            H \ar[r]^{\beta_A(x)}  \ar@{}[dr]|{\widetilde{\psi} \Downarrow}  \ar[d]_y 
          & A \ar[d]|{o}^{w_A}
          \\
            B  \ar[r]|{o}_{w_B}
          & E 
         }$. Since $\varphi \sim \widetilde{\varphi}$ and $\psi \sim \widetilde{\psi}$, we have the isomorphism
         
         $$
         ( 
\xymatrix@-.5pc@ur
         {
            G \ar[r]^x  \ar@{}[dr]|{id \Downarrow}  \ar[d]_{G(w_A)(x)}
          & A \ar[d]|{o}^{w_A}
          \\
            E  \ar[r]|{o}_{id}
          & E 
         },
\xymatrix@-.5pc@ur
         {
            H \ar[r]^y  \ar@{}[dr]|{id \Downarrow}  \ar[d]_{H(w_B)(y)}
          & A \ar[d]|{o}^{w_A}
          \\
            E  \ar[r]|{o}_{id}
          & E 
         } 
         \;) : (x,y,\varphi,\psi) \rightarrow (G(w_A)(x),H(w_B)(y),\widetilde{\varphi},\widetilde{\psi}). $$
         
         But $(G(w_A)(x),H(w_B)(y),\widetilde{\varphi},\widetilde{\psi}) = \Diamond (G(w_A)(x),H(w_B)(y),\widetilde{\varphi},\widetilde{\psi})$.

\vspace{1.5ex}         
         
\noindent \textit{(b) faithful.}  It is immediate from Proposition \ref{prop:mismoalphaybeta}.

\vspace{1.5ex}

\noindent \textit{(c) full.} Consider a premorphism of $Eq$ between two objects in the image of $\Diamond$, that is 

$$
( \xymatrix@ur@-.5pc{G \ar[r]^x \ar[d]_{x'} \ar@{}[dr] |{\Downarrow \xi} & A \ar[d]|{o}^{u} \\ A' \ar[r]_{v} & C},
\xymatrix@ur@-.5pc{H \ar[r]^y \ar[d]_{y'} \ar@{}[dr] |{\Downarrow \eta} & A \ar[d]|{o}^r \\ A' \ar[r]_{s} & D} 
\,): (x,y,\lambda_A(\gamma),\lambda_A(\delta)) \mr{} (x',y',\lambda_{A'}(\gamma'),\lambda_{A'}(\delta')).$$
Then the following two equations hold in $\cc{L}_{\sigma}(H)$
 \begin{align*} \eta \circ \lambda_A(\gamma) \sim \lambda_{A'}(\gamma') \circ \alpha_C(\xi)\\
 \eta \circ \lambda_A(\delta) \sim \lambda_{A'}(\delta') \circ \beta_C(\xi)
\end{align*}
Applying Lemma \ref{lema:todofinitov2}, we have $X \mr{w_X} E$ for $X=A,A',C,D$, 
and equivalent premorphisms such that the following equations hold in $HE$         
 \begin{align} \label{align:equationsabove} \widetilde{\eta} \circ \widetilde{\lambda_A(\gamma)} = \widetilde{\lambda_{A'}(\gamma')} \circ \widetilde{\alpha_C(\xi)} \\
 \label{align:equationsabove2} \widetilde{\eta} \circ \widetilde{\lambda_A(\delta)} = \widetilde{\lambda_{A'}(\delta')} \circ \widetilde{\beta_C(\xi)}
\end{align}
As in the proof of Theorem \ref{teo:conmutanconcotensors}, item (c), by \ref{lema:todofinitov2}, item 2 i), we have that $\widetilde{\lambda_A(\gamma)} = H(w_A)(\gamma)$, and similarly for the other three possibilities that appear in the equations above involving $A'$ and $\delta$.

Also, by \ref{lema:todofinitov2}, item 2 ii), $\widetilde{\alpha_C(\xi)}$ and $\widetilde{\beta_C(\xi)}$ are computed from $\alpha_C(\xi)$ and $\beta_C(\xi)$ pasting the same $2$-cells $\mu, \nu$ as in \eqref{eq:premorfismostilde}. Then, if we construct $\widetilde{\xi}$ from $\xi$ by pasting these $2$-cells $\mu, \nu$, since $\alpha$ and $\beta$ are $2$-natural transformations, we have $\widetilde{\alpha_C(\xi)} = \alpha_E(\widetilde{\xi})$, $\widetilde{\beta_C(\xi)} = \beta_E(\widetilde{\xi})$.

Then the equations \eqref{align:equationsabove}, \eqref{align:equationsabove2} above become
 \begin{align*} \widetilde{\eta} \circ H(w_A)(\gamma) &= H(w_{A'})(\gamma') \circ \alpha_E(\widetilde{\xi})\\
 \widetilde{\eta} \circ H(w_A)(\delta) &= H(w_{A'})(\delta') \circ \beta_E(\widetilde{\xi})
\end{align*}
 These equations show that $(\tilde{\xi},\tilde{\eta})$ is a premorphism in $\cc{L}_{\sigma}(P)$ between $(x,y,\gamma,\delta)$ and $(x',y',\gamma',\delta')$
\end{proof}


This finishes the proof of Theorem \ref{teo:conmutan}.

\bibliographystyle{unsrt}

\end{document}